\documentclass[onefignum,onetabnum]{siamart171218}
\usepackage{lipsum}
\usepackage{amsfonts}
\usepackage{graphicx}
\usepackage{epstopdf}
\usepackage{algorithmic}
\usepackage{siunitx} 
\usepackage{multirow}
\usepackage{tikz}
\usepackage{threeparttable}
\usepackage{amssymb}
\usepackage{caption}
\usepackage{subcaption}
\ifpdf
  \DeclareGraphicsExtensions{.eps,.pdf,.png,.jpg}
\else
  \DeclareGraphicsExtensions{.eps}
\fi

\newsiamremark{remark}{Remark}
\newsiamremark{hypothesis}{Hypothesis}
\crefname{hypothesis}{Hypothesis}{Hypotheses}
\newsiamthm{claim}{Claim}

\headers{Matrix-free parallel deflation method for Helmholtz problems}{J. Chen, V. Dwarka, and C. Vuik}

\title{A Matrix-free parallel two-level deflation preconditioner for the two-dimensional Helmholtz problems}

\author{Jinqiang Chen\thanks{Delft Institute of Applied Mathematics, TU Delft, the Netherlands 
  (\email{j.chen-11@tudelft.nl}, \email{v.n.s.r.dwarka@tudelft.nl}, \email{c.vuik@tudelft.nl}).}
\and Vandana Dwarka\footnotemark[1]
\and Cornelis Vuik\footnotemark[1]}

\usepackage{amsopn}

\makeatletter
\newcommand*{\addFileDependency}[1]{% argument=file name and extension
  \typeout{(#1)}% latexmk will find this if $recorder=0 (however, in that case, it will ignore #1 if it is a .aux or .pdf file etc and it exists! if it doesn't exist, it will appear in the list of dependents regardless)
  \@addtofilelist{#1}% if you want it to appear in \listfiles, not really necessary and latexmk doesn't use this
  \IfFileExists{#1}{}{\typeout{No file #1.}}% latexmk will find this message if #1 doesn't exist (yet)
}
\makeatother

\newcommand{\fineu}{ u^{h} }
\newcommand{\coarseu}{ u^{2h} }
\ifpdf
\hypersetup{
  pdftitle={A Matrix-free parallel two-level-deflation preconditioner for the two-dimensional Helmholtz problems},
  pdfauthor={J. Chen, V. Dwarka, and C. Vuik}
}
\fi

\begin{document}
\maketitle

% REQUIRED
\begin{abstract}
  We propose a matrix-free parallel two-level-deflation preconditioner combined with the Complex Shifted Laplacian preconditioner(CSLP) for the two-dimensional Helmholtz problems. The Helmholtz equation is widely studied in seismic exploration, antennas, and medical imaging. It is one of the hardest problems to solve both in terms of accuracy and convergence, due to scalability issues of the numerical solvers. Motivated by the observation that for large wavenumbers, the eigenvalues of the CSLP-preconditioned system shift towards zero, deflation with multigrid vectors, and further high-order vectors were incorporated to obtain wave-number-independent convergence. For large-scale applications, high-performance parallel scalable methods are also indispensable. In our method, we consider the preconditioned Krylov subspace methods for solving the linear system obtained from finite-difference discretization. The CSLP preconditioner is approximated by one parallel geometric multigrid V-cycle. For the two-level deflation, the matrix-free Galerkin coarsening as well as high-order re-discretization approaches on the coarse grid are studied. The results of matrix-vector multiplications in Krylov subspace methods and the interpolation/restriction operators are implemented based on the finite-difference grids without constructing any coefficient matrix. These adjustments lead to direct improvements in terms of memory consumption. Numerical experiments of model problems show that wavenumber independence has been obtained for medium wavenumbers. The matrix-free parallel framework shows satisfactory weak and strong parallel scalability.

 %1486-character version
 %We propose a matrix-free parallel two-level-deflation preconditioner combined with the Complex Shifted Laplacian preconditioner(CSLP) for the two-dimensional Helmholtz problems. The Helmholtz problem, which is widely studied in seismic exploration, is hard to solve both in terms of accuracy and convergence, due to the scalability issues of the numerical solvers. Motivated by the observation that for large wavenumbers, the eigenvalues of the CSLP-preconditioned system shift towards zero, deflation with multigrid vectors and further high-order vectors were incorporated to obtain wave-number-independent convergence. For large-scale applications, high-performance parallel scalable methods are also indispensable. In our method, we use the preconditioned Krylov subspace methods to solve the linear system obtained from finite-difference discretization. The CSLP preconditioner is approximately inverted by one parallel geometric multigrid V-cycle. Instead of the Galerkin coarsening method, we choose re-discretization on the coarse grid. The matrix-vector products and the interpolation/restriction operations are implemented based on the finite-difference grids without constructing the coefficient matrix. These adjustments lead to direct improvements in terms of memory consumption. Numerical experiments show that wavenumber independence has been obtained for medium wavenumbers. The matrix-free parallel framework shows satisfactory parallel performance and weak scalability.

\end{abstract}

% REQUIRED

% REQUIRED, I don't know what is it now

\section{Introduction}
The Helmholtz equation, describing the phenomena of time-harmonic wave scattering in the frequency domain, finds applications in many scientific fields, such as seismic problems, advanced sonar devices, medical imaging, and many more. To solve the Helmholtz equation numerically, we discretize it and obtain a linear system $Ax=b$. The matrix of the linear system is sparse, symmetric, complex-valued, non-Hermitian, and indefinite. Instead of a direct solver, iterative methods and parallel computing are commonly considered for a large-scale linear system resulting from a practical problem. However, the indefiniteness of the linear system brings a great challenge to the numerical solution method, especially for large wavenumbers. The convergence rate of many iterative solvers is affected significantly by increasing wavenumber. An increase in the wavenumber leads to a dramatic increase in iterations. Moreover, the general remedy for minimizing the so-called pollution error, driven by numerical dispersion errors due to discrepancies between the exact and numerical wavenumber, is to refine the grid such that the condition $k^3h^2 < 1$ is satisfied \cite{babuska1997pollution}. Therefore, in this field, the research problem of how to solve the systems efficiently and economically while at the same time maintaining high accuracy by minimizing pollution errors arises. A wavenumber-independent-convergent and parallel scalable iterative method could significantly enhance the corresponding research in electromagnetic, seismology, and acoustics.

Many efforts have been made to solve the problem in terms of accuracy and scalable convergence behavior. One of the main concerns is the spectrum of the system matrix, which is closely related to the convergence of Krylov subspace methods. The preferable idea is to preprocess the system with a preconditioner. By applying a preconditioner to the linear system, the solution remains the same, but the coefficient matrix has a more favorable distribution of eigenvalues. 
Many preconditioners are proposed for the Helmholtz problem so far, such as incomplete Cholesky (IC)/incomplete LU (ILU) factorization, shifted Laplacian preconditioners \cite{bayliss1983iterative,erlangga2004class,erlangga2006novel}, and so on. The Complex Shifted Laplace Preconditioner (CSLP) \cite{erlangga2004class, erlangga2006novel} does show good properties for medium wavenumbers. Nevertheless, the eigenvalues shift to the origin as the wavenumber increases. The convergence speed of Krylov-based iterative solvers is adversely affected by the presence of near-zero eigenvalues. Although spectral analysis for a nonnormal matrix becomes less meaningful to assess convergence properties, numerical experiments show that this notion is also true. \cite{dwarka2020scalable}.

The deflation method was developed to eliminate unfavorable eigenvalues by projecting them to zero. Initially proposed independently by Nicolaides \cite{nicolaides1987deflation} and Dost{'a}l \cite{dostal1988conjugate} as a means to accelerate the standard conjugate gradient method for symmetric positive definite (SPD) systems, the deflation method has been further investigated by \cite{mansfield1991damped, kolotilina1998twofold, saad2000deflated, tang2008two}. The deflation preconditioner is sensitive to approximations of the coarse-grid system, but obtaining the exact inverse of the coarse-grid matrix for large-scale problems is not always feasible. In the absence of an accurate approximation, the projection step can introduce numerous new eigenvalues close to zero, leading to adverse effects. To address this limitation, Erlangga et al. \cite{erlangga2008multilevel} extended the method to handle systems with nonsymmetric matrices. Instead of projecting the small eigenvalues to zero, they deflated the smallest eigenvalues to match the maximum eigenvalue. Erlangga et al. demonstrated that their modified projection method was less sensitive to approximations of the coarse-grid system, allowing for the utilization of multilevel projected Krylov subspace iterations. In a recent development, Dwarka and Vuik \cite{dwarka2020scalable} introduced higher-order approximation schemes for constructing deflation vectors. This advanced two-level deflation method exhibits convergence that is nearly independent of the wavenumber. The authors further extend the two-level deflation method to a multilevel deflation method \cite{DWARKA2022111327}. By using higher-order deflation vectors, they show that up to the level where the coarse-grid linear systems remain indefinite, the near-zero eigenvalues of these coarse-grid operators remain aligned with the near-zero eigenvalues of the fine-grid operator. This keeps the spectrum of the preconditioned system away from the origin. Combining this with the well-known CSLP preconditioner, they obtain a scalable solver for highly indefinite linear systems.

The development of scalable parallel Helmholtz solvers is also ongoing. One approach is to parallelize existing advanced algorithms. Kononov and Riyanti \cite{Kononov2007Numerical,riyanti2007parallel} first developed a parallel version of Bi-CGSTAB preconditioned by multigrid-based CSLP. Dan and Rachel \cite{Gordon2013Robust} parallelized their so-called CARP-CG algorithm (Conjugate Gradient Acceleration of CARP) blockwise. The block-parallel CARP-CG algorithm shows improved scalability as the wavenumber increases. Calandra et al. \cite{calandra2013improved,calandra2017geometric} proposed a geometric two-grid preconditioner for 3D Helmholtz problems, which shows the strong scaling property in a massively parallel setup.

Another approach is the Domain Decomposition Method (DDM), which originates from the early Schwarz methods. DDM has been widely used to develop parallel solution methods for Helmholtz problems. For comprehensive surveys, we refer the reader to \cite{douglas1998second,mcinnes1998additive,toselli1999overlapping,collino2000domain,gander2002optimized,schadle2007additive,engquist2011sweeping,boubendir2012quasi,chen2013source,stolk2013rapidly,gander2019class,taus2020sweeps} and references therein.

%%\lipsum[2-3]
Based on the parallel framework of the CSLP-preconditioned Krylov subspace methods \cite{jchen2D2022}, this paper proposes a matrix-free parallel implementation of the two-level-deflation preconditioner for the two-dimensional Helmholtz problems. In addition to the matrix-free parallelization of each component of the deflation technique, as one of our main contributions, we extensively studied and compared several matrix-free implementations of coarse-grid operators. We show that the coarse-grid re-discretization schemes derived from the Galerkin coarsening approach yield convergence properties that are nearly independent of the wavenumber. Furthermore, we investigate the weak and strong scaling properties of the present parallel solution method in a massively parallel setting.

% The outline is not required, but we show an example here.
The paper is organized as follows. The mathematical models and numerical methods are given in \cref{sec:MPs} and \cref{sec:NMs}. Our parallel implementation is in \cref{sec:imp}. \Cref{sec:matrix-free} introduces the matrix-free operators in detail. The experimental results will be discussed in \cref{sec:experiments}, and the conclusions follow in \cref{sec:conclusions}.

\section{Mathematical Models}
\label{sec:MPs}

In this paper, we mainly consider the following Helmholtz equation. Suppose that the domain $\Omega$ is rectangular with boundary $\Gamma = \partial\Omega$. The Helmholtz equation reads
	
\begin{equation}\label{eq:HelmholtzEq}
-\Delta \mathbf{u} - k(x,y)^2 \mathbf{u} = \mathbf{b},\  \text{on} \ \Omega 
\end{equation}
supplied with one of the following boundary conditions:
\begin{equation}\label{eq:DirichletBC}
\text{Dirichlet: } \mathbf{u}=\mathbf{g},\  \text{on} \ \partial\Omega 
\end{equation}

\begin{equation}\label{eq:SommerfeldBC}
\text{first-order Sommerfeld: } \frac{\partial \mathbf{u}}{\partial \vec{n}}-\text{i} k(x,y) \mathbf{u} = \mathbf{0},\  \text{on} \ \partial\Omega 
\end{equation}
where $\text{i}$ is the imaginary unit. $\vec{n}$ and $\mathbf{g}$ represent the outward normal and the given data of the boundary, respectively. $\mathbf{b}$ is the source function. $k(x,y)$ is the wavenumber on $\Omega$. Suppose the frequency is $f$, the speed of propagation is $c(x,y)$, they are related by
\begin{equation}\label{eq:wavenum}
k(x,y) =\frac{2\pi f}{c(x,y)}.
\end{equation}

\subsection{Discretization} \label{sec:Discretization}
Structural vertex-centered grids are used to discretize the computational domain. Suppose the mesh width in $x$ and $y$ direction are both $h$.  A second-order finite difference scheme for a 2D Laplace operator has the following stencil
\begin{equation}\label{eq:stencilLaplace}
\left[ {{-\Delta_h}} \right] = \frac{1}{{{h^2}}}\left[ {\begin{array}{*{20}{c}}
    0&{ - 1}&0\\
    { - 1}&{4}&{ - 1}\\
    0&{ - 1}&0
    \end{array}} \right].
\end{equation}
The discrete Helmholtz operator $A_h$ can be obtained by subtracting the diagonal matrix $ \mathcal{I}(k_{i,j}^2)_h$  to the Laplacian operator $-\Delta_h$, \textit{i.e.}
\begin{equation}\label{eq:operator}
A_h =-\Delta_h- \mathcal{I}(k_{i,j}^2)_h
\end{equation}
Therefore, the stencil of the discrete Helmholtz operator is
\begin{equation}\label{eq:Helmstencil}
\left[ {{A_h}} \right] = \frac{1}{{{h^2}}}\left[ {\begin{array}{*{20}{c}}
    0&{ - 1}&0\\
    { - 1}&{4 - {k_{i,j}^2}{h^2}}&{ - 1}\\
    0&{ - 1}&0
    \end{array}} \right]
\end{equation}

In case the Sommerfeld radiation condition \eqref{eq:SommerfeldBC} is used, the discrete schemes for the boundary points need to be defined.  Ghost points located outside the boundary points can be introduced. For instance, suppose $u_{0,j}$ is a ghost point on the left of $u_{1,j}$, the normal derivative can be approximated  by
\begin{equation}
    \label{eq:GhostPoint}
    \frac{\partial u}{\partial \vec{n}} -\text{i} k(x,y) u =\frac{u_{0,j}-u_{2,j}}{2h}-\text{i} k_{1,j} u_{1,j}= 0
\end{equation}
We can rewrite it as
\begin{equation}
u_{0,j} = u_{2,j} + 2h\text{i} k_{1,j} u_{1,j}.
\end{equation}
With the elimination of the ghost point, the stencil for the boundary points $u_{1,j}$ becomes
\begin{equation}\label{eq:Sommerfeldstencil1}
    \frac{1}{{{h^2}}}\left[ {\begin{array}{*{20}{c}}
    0&{ - 1}&0\\
    { 0}&{4 - {k_{1,j}^2}{h^2}-2\text{i} k_{1,j}h}&{ - 2}\\
    0&{ - 1}&0
    \end{array}} \right]
\end{equation}
For corner points, we can do the above step in both $x$- and $y$- directions. The resulting stencil is
\begin{equation}\label{eq:Sommerfeldstencil2}
\frac{1}{{{h^2}}}\left[ {\begin{array}{*{20}{c}}
    0&{ - 2}&0\\
    { 0}&{4 - {k_{1,1}^2}{h^2}-4\text{i} k_{1,1}h}&{ - 2}\\
    0&{ 0}&0
    \end{array}} \right]
\end{equation}

Discretization of the partial equation on the finite-difference grids results in a system of linear equations 
\begin{equation}\label{eq:linearsys}
A_h \mathbf{u}_h=\mathbf{b}_h
\end{equation}
With first-order Sommerfeld boundary conditions, the resulting matrix is sparse, symmetric, complex-valued, indefinite, and non-Hermitian. % for a sufficiently large wavenumber $k$.

Note that $kh$ is an important parameter that indicates how many grid points per wavelength are needed. The mesh width  $h$ can be determined by the guidepost of including at least $N_{pw}$(e.g. 10 or 30) grid points per wavelength. They have the following relationships 
\begin{equation}
kh=\frac{2 \pi h}{\lambda} = \frac{2 \pi }{N_{pw}}.
\end{equation}
For example, if at least 10 grid points per wavelength are required, we can maintain $kh = 0.625$.

% \subsection{Model Problem 1}
% This problem is the so-called 2D closed-off problem. A rectangular homogeneous domain $\Omega=\left[0,1 \right]$ is considered. The source function is specified by
% \begin{equation}\label{eq:2Dcloseoff}
% b\left( {x,y} \right) = \left( {5{\pi ^2} - {k^2}} \right)\sin \left( {\pi x} \right)\sin \left( {2\pi y} \right) - {k^2},\ \Omega=\left[0,1 \right]  
% \end{equation}
% It is supplied with the following Dirichlet conditions
% \begin{equation}
% u=1,\  \text{on} \ \partial\Omega \nonumber
% \end{equation}

% Thus, we can have the exact solution given by
% \begin{equation}
% u(x,y) = \sin \left( {\pi x} \right)\sin \left( {2\pi y} \right) + 1.
% \end{equation}

\subsection{Model Problem - constant wavenumber}
We first consider a 2D problem with constant wavenumber in a rectangular homogeneous domain $\Omega=\left[0,1 \right]$. A point source is given by  
\begin{equation}\label{eq: dirac rhs}
b\left( {x,y} \right) = \delta \left( {x - x_0,y - y_0} \right),\ \Omega=\left[0,1 \right]  
\end{equation}
where $\delta \left( {x,y} \right)$ is a Dirac delta function in the following form
\begin{equation}\label{eq: dirac func}
\delta \left( {x,y} \right) = \left\{ {\begin{array}{*{20}{c}}
    { + \infty \;\;\;x = 0,\;y = 0}\\
    {0\;\;\;\;\;\;\;x \ne 0,\;y \ne 0}
    \end{array}} \right.
\end{equation}
which satisfies
\begin{equation}\label{eq: dirac func restricton}
\int {\int {{\delta ^2}\left( {x,y} \right)dxdy = 1} } 
\end{equation}

In this problem, the point source is imposed at the center $(x_0,y_0)=(0.5,0.5)$. The wave propagates outward from the center of the domain. The Dirichlet boundary conditions (denoted as MP-2a) or the first-order Sommerfeld radiation conditions (denoted as MP-2b) are imposed, respectively.  

\subsection{Model Problem - Wedge problem}
\begin{figure}[htbp]
    \centering
    \includegraphics[width=0.38\textwidth]{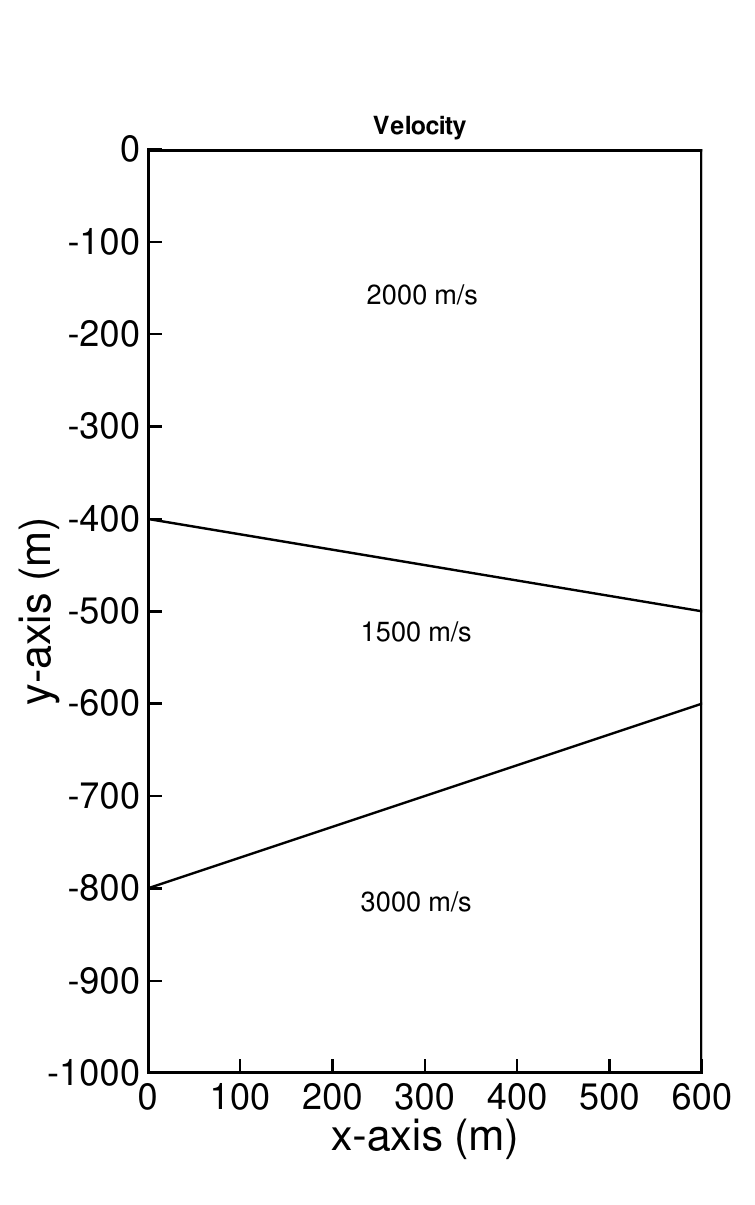}
    \caption{The velocity distribution of the wedge problem}
    \label{fig:wedge domain}
\end{figure}

Most physical problems of geophysical seismic imaging describe a heterogeneous medium. The so-called Wedge problem \cite{plessix2003separation} is a typical problem with a simple heterogeneity. It mimics three layers with different velocities hence different wavenumbers. As shown in Figure \ref{fig:wedge domain}, the rectangular domain $\Omega=\left[0,600 \right] \times \left[-1000,0 \right]$  is split into three layers. Suppose the wave velocity $c$ is constant within each layer but different from each other. A point source is located at $\left(x, y \right) = \left( 300, 0\right) $. 
	
The problem is given by
\begin{equation}\label{eq: wedge problem}
    \left\{ {\begin{array}{*{20}{c}}
        {-\Delta u(x,y)-k(x,y)^2u(x,y)=b(x,y), \quad\text{on } \Omega=(0,600)\times(-1000,0)}\\
        {b(x,y)=\delta(x-300,y) \quad x,y\in \Omega}
        \end{array}} \right.
\end{equation}
where $k(x,y)=\frac{2\pi f}{c(x,y)}$. The wave velocity $c(x,y)$ is shown in Figure \ref{fig:wedge domain}. $f$ is the frequency. The first-order Sommerfeld boundary conditions are imposed on all boundaries.

\subsection{Model Problem - Marmousi Problem} \label{sec:Marmousi}
	For industrial applications, the fourth model problem is the so-called Marmousi problem \cite{Versteeg_1991_ME}, a well-known benchmark problem. The geometry of the problem stems from the section of the Kwanza Basin through North Kungra. It contains 158 horizontal layers in the depth direction, making it highly heterogeneous. 

\begin{figure}[htbp]
		\centering
		\includegraphics[width=0.89\textwidth, trim=4 4 4 4,clip]{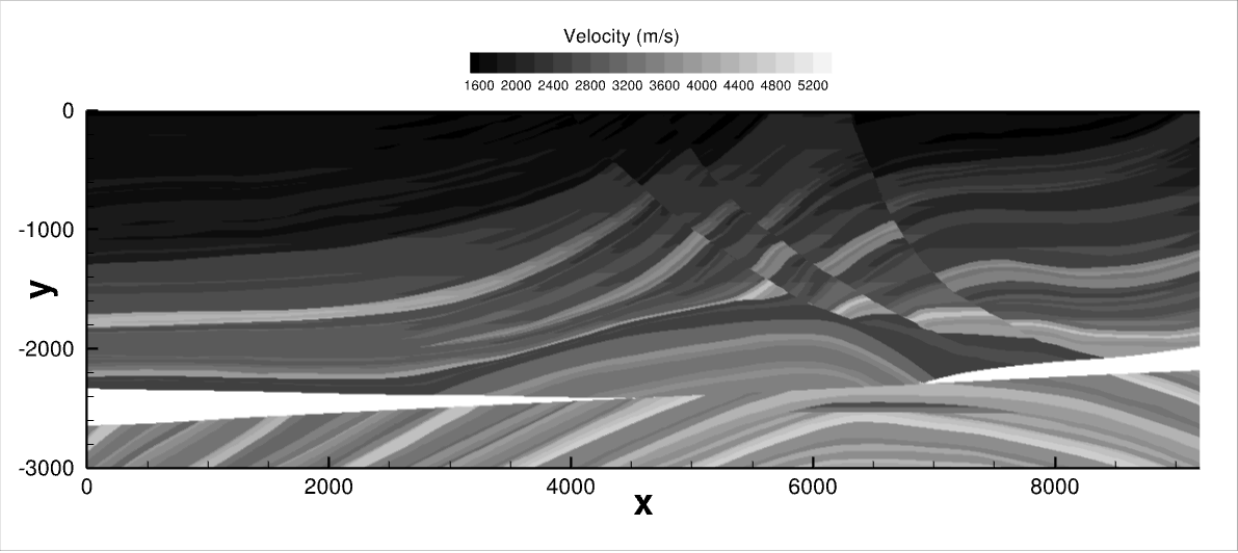}
		\caption{The velocity distribution  of the marmousi problem}
		\label{fig:marmousi}
	\end{figure}

	The problem is given by
	\begin{equation}\label{eq: marmousi problem}
		\left\{ {\begin{array}{*{20}{c}}
			{-\Delta u(x,y)-k(x,y)^2u(x,y)=b(x,y), \quad\text{on } \Omega=(0,9200)\times(-3000,0)}\\
			{b(x,y)=\delta(x-6000,y) \quad x,y\in \Omega}
			\end{array}} \right.
	\end{equation}
	where $k(x,y)=\frac{2\pi f}{c(x,y)}$. The wave velocity $c(x,y)$ over the domain is shown in Figure \ref{fig:marmousi}. The first-order Sommerfeld boundary conditions are imposed on all boundaries.

\section{Numerical Methods}
\label{sec:NMs}
The Krylov subspace methods are popular for solving large and sparse linear systems. The coefficient matrix of a non-Hermitian linear system is preferred to have a spectrum located in a bounded region that excludes the origin in the complex plane. Because the iterative methods can converge fast in such cases. We can incorporate preconditioning and deflation to enhance the convergence of Krylov subspace methods. This section will specify every component of the preconditioned Krylov subspace methods we used to solve the linear system  resulting from the finite-difference discretization. The multigrid-based Complex shifted Laplace preconditioner (CSLP) combined with two-level deflation will be considered in this work.

\subsection{Complex Shifted Laplace Preconditioner}
We focus on the CSLP due to its superior performance and easy setup. 

The CSLP is defined by
\begin{equation}\label{CSLP def}
    M_{h,(\beta_1,\beta_2)}=-\Delta_h-(\beta_1+\text{i}\beta_2) \mathcal{I}(k_{i,j}^2)_h
\end{equation}
where $\text{i}=\sqrt{-1}$ and $(\beta_1,\beta_2)\in [0,1]$. The standard multigrid methods \cite{maclachlan2008algebraic,erlangga2006novel} are usually employed to invert a shifted Laplacian preconditioner. In the numerical experiments of this report, unless noted otherwise, $\beta_1=1, \beta_2=-0.5$ will be used \cite{Gijzen2007,gander2015applying}. The inverse of the preconditioner is approximated by a multigrid V-cycle, which consists of one step of pre-/post- damped Jacobi smoother with a relaxation parameter of $0.8$, full weighting restriction, bilinear interpolation, and a General Minimal RESidual method (GMRES) \cite{saad1986gmres} solver for solving the coarsest-grid problem.

\subsection{Deflation}
Suppose a general linear system $Au=b$, where $A\in\mathbb{R}^{n \times n}$, and a projection subspace matrix, $Z \in \mathbb{R}^{n\times m}$, with $m<n$ and full rank are given. Assume that $E=Z^T A Z$ is invertible, the projection matrix $P$ can be defined as
\begin{equation}
    P = I - AQ, \quad Q=Z E^{-1} Z^T, \quad E = Z^T A Z 
\end{equation}
By applying the deflation operator, the linear system to be solved becomes
\begin{equation}
\label{eq: PAu eq Pb}
    PA\hat{u}=Pb
\end{equation}
Krylov subspace method can be employed to solve the resulting linear system. Since equation (\ref{eq: PAu eq Pb}) is singular, to obtain a unique solution, we need to update $u$ by
\begin{equation}
    u = Qb+P^T\hat{u}
\end{equation}

\subsubsection{Deflation Vectors for Helmholtz}
The choice of the deflation vectors is one of the base components of the deflation preconditioner. For the Helmholtz problem, a good and necessarily sparse approximation of the eigenvectors corresponding to the small eigenvalues is ideal. However, it is computationally expensive. Originating from the observation that the multigrid inter-grid operators highlight the small frequencies and reserve them on the coarser level, we employ the geometrically constructed multigrid vectors as the deflation vectors. The standard linear interpolation along one dimension is given by 
\begin{equation}
\label{eq: linear iterpolation}
I_{2h}^{h} \coarseu_i=\left\{\begin{array}{cl}
\coarseu_{i/2} & \text{if } i\text{ is even}\\
\frac{1}{2}\left( \coarseu_{\left(i-1\right)/2} + \coarseu_{\left(i+1\right)/2}\right) & \text{if } i\text{ is odd}
\end{array}\right.
\end{equation}
for $i=1,2,\dots,n-1$ and the restriction operator is
\begin{equation}
\label{eq: linear restriction}
I_{h}^{2h} \fineu_i=\frac{1}{2}\left( \fineu_{2i-1} +2\fineu_{2i}+\fineu_{2i+1} \right) 
\end{equation}
for $i=1,2,\dots,n/2$. 
In two-dimension, the so-called full weighting restriction operator has a stencil given by
\begin{equation}\label{eq:weight operator}
[I_h^{2h}] = \frac{1}{16}\left[ {\begin{array}{*{20}{c}}
    1&2&1\\
    2&4&2\\
    1&2&1
    \end{array}} \right]_h^{2h},
\end{equation}
and the stencil of the bilinear interpolation operator is given by
\begin{equation}\label{eq:bilinear operator}
[I_{2h}^h] = \frac{1}{4}\left[ {\begin{array}{*{20}{c}}
    1&2&1\\
    2&4&2\\
    1&2&1
    \end{array}} \right]_{2h}^h.
\end{equation}
% For a two-dimensional setup, we can perform the bilinear interpolation as the following polynomials
% \begin{equation}
% \label{eq: bilinear iterpolation}
%     \left[ I_{2h}^{h}\right] = \frac{1}{4}\begin{bmatrix}
%         1 & 2 & 1\\
%         2 & 4 & 2\\
%         1 & 2 & 1
%     \end{bmatrix}_{2h} ^{h}
% \end{equation}
If we set the coarse to fine grid interpolation operator as the deflation subspace
\begin{equation}
    Z = I_{2h}^{h}
\end{equation}
where $h \rightarrow 2h$ indicates the standard coarsening method. Then $Z^T$ will be the full-weighting restriction operator. The deflation preconditioner can be rewritten as
\begin{equation}
    P_h = I_h - A_h Q_h, \quad Q_h=I_{2h}^{h} A_{2h}^{-1} I_{h}^{2h}, \quad A_{2h} = I_{h}^{2h} A_{h} I_{2h}^{h} 
\end{equation}

Recently, Dwarka and Vuik \cite{dwarka2020scalable} used higher-order B\'ezier curves to construct deflation vectors and obtained wave-number-independent convergence for very large wavenumbers. The 1D interpolation and restriction operator is given by
\begin{subequations}
	\begin{equation}
	I_{2h}^{h} \coarseu_i=\begin{cases}
	\frac{1}{8}\left( \coarseu_{\left(i-2\right)/2} +6 \coarseu_{i/2}+\coarseu_{\left( i+2\right)/2}\right) & \text{if } i\text{ is even}\\
	\frac{1}{2}\left( \coarseu_{\left(i-1\right)/2} + \coarseu_{\left(i+1\right)/2}\right) & \text{if } i\text{ is odd}
	\end{cases}
	\label{eq:quadratic approximation-a}
	\end{equation}
	\begin{equation}
	I_{h}^{2h} \fineu_i=\frac{1}{8}\left( \fineu_{ 2i-2}+4\fineu_{2i-1} + 6\fineu_{2i}+4\fineu_{2i+1}+\fineu_{ 2i+2 } \right)
	\label{eq:quadratic approximation-b}
	\end{equation}
	\label{eq:quadratic approximation}
\end{subequations}

\begin{figure}[htb]
  \centering
  \begin{tikzpicture}[scale=1.0]
	
    \foreach \x in {0,1,2,3,4}
    \foreach \y in {0,1,2,3,4}
    {
    %#####################################################
    \ifthenelse{  \lengthtest{\x pt < 4pt}  }{
    \draw[dashed] (\x,\y) -- (\x+1,\y);
    }{}
    %#####################################################
    \ifthenelse{  \lengthtest{\y pt < 4pt}  }{
    \draw[dashed] (\x,\y) -- (\x,\y+1);
    }{}	
    }

    \foreach \x in {0,2,4}
    \foreach \y in {0,2,4}
    {
    %#####################################################
    \ifthenelse{  \lengthtest{\x pt < 4pt}  }{
    \draw (\x,\y) -- (\x+2,\y);
    }{}
    %#####################################################
    \ifthenelse{  \lengthtest{\y pt < 4pt}  }{
    \draw (\x,\y) -- (\x,\y+2);
    }{}
    }
        \foreach \x in {0,2,4}
        \foreach \y in {0,2,4}
       {
         \filldraw[color=black, fill= black] (\x,\y) circle (1.5pt);
        }

        \foreach \x in {1,3}
        \foreach \y in {0,2,4}
       {
        \filldraw[dashed, color=black, fill= blue] (\x,\y) circle (1.5pt);
        }

        \foreach \y in {1,3}
        \foreach \x in {0,2,4}
       {
        \filldraw[dashed, color=black, fill= red] (\x,\y) circle (1.5pt);
        }

        \foreach \y in {1,3}
        \foreach \x in {1,3}
       {
        \filldraw[dashed, color=black, fill= black!20] (\x,\y) circle (1.5pt);
        }

        \filldraw[color=black, fill= black] (0.0,-0.8) circle (1.5pt);
        \filldraw[dashed,color=black, fill= blue] (0.2,-0.8) circle (1.5pt);
        \filldraw[dashed,color=black, fill= red] (0.4,-0.8) circle (1.5pt);
        \filldraw[dashed,color=black, fill= black!20] (0.6,-0.8) circle (1.5pt);
        \draw (0.65,-0.8) node[scale=1, right] {: fine grid points $\in\Omega^{h}$};
        
        \filldraw[color=black, fill= black] (0.3,-1.2) circle (1.5pt);			
        % \shade[blueBall] (1,-1) circle (0.05cm);
      %\node at (3,0,5) [circle,fill=black] {};
      \draw (0.65,-1.2) node[scale=1, right] {: coarse grid points $\in\Omega^{2h}$};
      \draw (2.01,1.99) -- (2.08,1.85);
      \draw (2.05,1.75) node[scale=0.5, right] {$(i,j)\in \Omega^{2h}$};
      %\draw (2.05,1.75) node[scale=0.5, right] {$(2i-1,2j-1,2k-1)\in \Omega^{h}$};

  \end{tikzpicture}
  \caption{The allocation map of interpolation operator}
  \label{fig: high order interpolation stencil}
\end{figure}
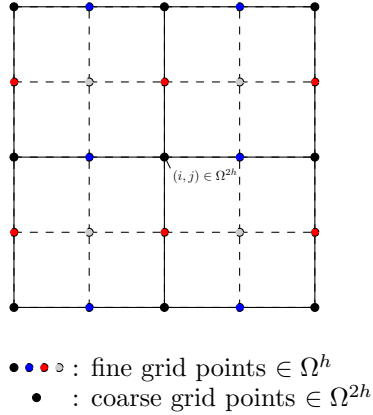

For two-dimensional cases, as shown in Figure \ref{fig: high order interpolation stencil}, the higher-order interpolation can be implemented as follows.
\begin{equation}\label{eq: O4 matrix-free prolongation}
    \begin{aligned}
        &I_{2h}^h u_{i_c,j_c}^{2h} = \\
        &\left\{ \begin{aligned}
        \frac{1}{64} & \left( 6 u_{i_c-1,j_c}^{2h} + 36 u_{i_c,j_c}^{2h} + 6 u_{i_c+1,j_c}^{2h}\right.\\
        & + u_{i_c-1,j_c-1}^{2h} + 6 u_{i_c,j_c-1}^{2h} + u_{i_c+1,j_c-1}^{2h}\\
        & + \left. u_{i_c-1,j_c+1}^{2h} + 6 u_{i_c,j_c+1}^{2h} + u_{i_c+1,j_c+1}^{2h}\right), &&\textcolor{black}{\bullet} \\
        \frac{1}{16} & \left( u_{i_c,j_c-1}^{2h} + 6 u_{i_c,j_c}^{2h} + u_{i_c,j_c+1}^{2h}\right.\\
        & + \left. u_{i_c+1,j_c-1}^{2h} + 6 u_{i_c+1,j_c}^{2h} + u_{i_c+1,j_c+1}^{2h}\right), &&\textcolor{blue}{\bullet} \\
         \frac{1}{16} & \left( u_{i_c-1,j_c}^{2h} + 6 u_{i_c,j_c}^{2h} + u_{i_c+1,j_c}^{2h}\right.\\
        & + \left. u_{i_c-1,j_c+1}^{2h} + 6 u_{i_c,j_c+1}^{2h} + u_{i_c+1,j_c+1}^{2h}\right), &&\textcolor{red}{\bullet} \\
        \frac{1}{4} & \left( u_{i_c,j_c}^{2h} + u_{i_c,j_c+1}^{2h}+u_{i_c+1,j_c}^{2h} + u_{i_c+1,j_c+1}^{2h} \right), &&\textcolor{black!20}{\bullet} 
        \end{aligned} \right.
    \end{aligned}
\end{equation}
where $(i_c,j_c) \in \Omega^{2h}$. The higher-order interpolation can then be represented by the following stencil notation
\begin{equation}\label{eq:higher-order intp-operator}
[I_{2h}^{h}] = \frac{1}{64}\left[ {\begin{array}{*{20}{c}}
    1&4&6&4&1\\
    4&16&24&16&4\\
    6&24&36&24&6\\
    4&16&24&16&4\\
    1&4&6&4&1
    \end{array}} \right]_{2h}^{h},
\end{equation}

The higher-order restriction can be implemented in a matrix-free way as follows.
	\begin{equation}\label{eq: O4 matrix-free restriction}
		\begin{aligned}
			& I_h^{2h} u_{2i_c-1,2j_c-1}^h = \\
            &\frac{1}{64} \left[ u_{2i_c-3,2j_c+1}^h + 4u_{2i_c-2,2j_c+1}^h + 6u_{2i_c-1,2j_c+1}^h + 4u_{2i_c,2j_c+1}^h + u_{2i_c+1,2j_c+1}^h \right.\\
			& + 4u_{2i_c-3,2j_c}^h + 16u_{2i_c-2,2j_c}^h + 24u_{2i_c-1,2j_c}^h + 16u_{2i_c,2j_c}^h + 4u_{2i_c+1,2j_c}^h \\
            & + 6u_{2i_c-3,2j_c-1}^h + 24u_{2i_c-2,2j_c-1}^h + 36u_{2i_c-1,2j_c-1}^h + 24u_{2i_c,2j_c-1}^h + 6u_{2i_c+1,2j_c-1}^h \\
            & + 4u_{2i_c-3,2j_c-2}^h + 16u_{2i_c-2,2j_c-2}^h + 24u_{2i_c-1,2j_c-2}^h + 16u_{2i_c,2j_c-2}^h + 4u_{2i_c+1,2j_c-2}^h \\
            & +  \left. u_{2i_c-3,2j_c-3}^h + 4u_{2i_c-2,2j_c-3}^h + 6u_{2i_c-1,2j_c-3}^h + 4u_{2i_c,2j_c-3}^h + u_{2i_c+1,2j_c-3}^h \right],
		\end{aligned}
	\end{equation}
where $(i_c,j_c) \in \Omega^{2h}$. The higher-order restriction can then be represented by the following stencil notation
\begin{equation}\label{eq:higher-order res-operator}
[I_h^{2h}] = \frac{1}{64}\left[ {\begin{array}{*{20}{c}}
    1&4&6&4&1\\
    4&16&24&16&4\\
    6&24&36&24&6\\
    4&16&24&16&4\\
    1&4&6&4&1
    \end{array}} \right]_h^{2h}.
\end{equation}

\subsubsection{The Deflation Preconditioner}
First, one may want to deflate the CSLP-preconditioned Helmholtz operator, as the eigenvalues of $M_{h,(\beta_1,\beta_2)}^{-1}A_h$ begin to shift to the origin as the wavenumber increases. These small eigenvalues can be projected to zero by the preconditioner given by
\begin{equation}
        \tilde{P}_h = I_h - \tilde{A}_h \tilde{Q}_h, \quad \tilde{Q}_h=I_{2h}^{h} \tilde{A}_{2h}^{-1} I_{h}^{2h}, \quad \tilde{A}_{2h}= I_{h}^{2h} \tilde{A}_{h} I_{2h}^{h}
\end{equation}
where $\tilde{A}_{h}= M_{h,(\beta_1,\beta_2)}^{-1}A_h$. One can notice that the coarse operator $\tilde{A}_{2h}$ needs to be inverted. Since the idea of deflation is to project the small eigenvalues to zero, the exact inversion of $\tilde{A}_{2h}$ is necessary. As the direct solver is not easy to parallelize and it is impractical for large problems, one has to turn to approximate solvers. However, an approximation of $\tilde{A}_{2h}$ will lead to a cluster of near-zero values in the preconditioned spectrum. Therefore, we choose to deflate toward the largest eigenvalues of the preconditioned system by adding a term. Thus, the so-called Two-Level Krylov method (TLKM) preconditioner \cite{PhDSheikh2014} reads as
\begin{equation}
    \tilde{P}_{h,\gamma} = \tilde{P}_h + \gamma \tilde{Q}_{h} = (I_h - \tilde{A}_h \tilde{Q}_h) + \gamma \tilde{Q}_{h}
\end{equation}
where $\gamma$ is usually set equal to the largest eigenvalues of the preconditioned system. Here we set $\gamma = 1$ since the CSLP preconditioner we will include leads to eigenvalues bounded by 1 in modulus. 
The preconditioned linear system to be solved is
\begin{equation}
    \tilde{P}_{h,\gamma} \tilde{A}_{h} u_h = \tilde{P}_{h,\gamma} \tilde{b}_{h}
\end{equation}
where $\tilde{b}_{h} = M_{h,(\beta_1,\beta_2)}^{-1}b_h$. Compared to the original Helmholtz system, the TLKM preconditioner is equivalent to 
\begin{equation}
    P_{h,TLKM} = \tilde{P}_{h,\gamma} M_{h,(\beta_1,\beta_2)}^{-1}
\end{equation}
One can observe that TLKM needs an application of $M_{h,(\beta_1,\beta_2)}^{-1}$ on the fine grid for every coarse grid iteration. This makes it too expensive to use.

Alternatively, one can deflate the Helmholtz operator. The preconditioner becomes
\begin{equation}
    P_{h,\gamma} = P_{h} + \gamma Q_{h} = (I_{h}-A_{h}Q_{h})+\gamma Q_h
\end{equation}
Combined with the standard preconditioner CSLP, which is commonly used in the Helmholtz system, a robust two-level preconditioned method, Adapted Deflation Variant 1 (A-DEF1) \cite{PhDTang2008} reads as 
\begin{equation}
    \label{eq: ADEF1}
    P_{h,A-DEF1} = M^{-1}_{h,(\beta_1,\beta_2)}P_h+Q_h
\end{equation}
The preconditioned linear system to be solved becomes
\begin{equation}
    \label{eq: ADEF linear sys}
    P_{h,A-DEF1}A_{h}u_{h} = P_{h,A-DEF1}b_h.
\end{equation}
Note that $P_{h,A-DEF1}A_{h}$ is non-singular, so Equation \ref{eq: ADEF linear sys} has a unique solution. When the higher-order deflation vector is employed, it will be denoted as Adapted Preconditined Deflation (APD) \cite{dwarka2020scalable}.

\subsection{Preconditioned Krylov subspace method}
For a large sparse linear system, the Krylov subspace methods are popular. 
%Some representative Krylov methods, like Conjugate Gradient(CG), CGNR, MINRES, BICG, Bi-CGSTAB, and GMRES, are developed so far. Among these, the CG method is a basic one. The error is minimized in the A-norm, and it only needs three vectors in memory during iterations. However, this algorithm is designed for a symmetric and positive definite system matrix. In contrast, 
Bi-CGSTAB and GMRES-type methods can be used for non-singular problems that are indefinite and non-symmetric as well. Compared with GMRES-type methods, Bi-CGSTAB has short recurrences and is easily parallelizable. But GMRES-type methods are more robust. Thus, GMRES-type methods are a suitable choice for the Helmholtz equation if the total number of iterations required is not too large.

The coefficient matrix of a non-Hermitian linear system is preferred to have a spectrum located in a bounded region that excludes the origin in the complex plane. Iterative methods can converge fast in such cases. We can incorporate a preconditioner to enhance the convergence of Krylov subspace methods. That is to pre-multiply the linear system with a preconditioning matrix $M^{-1}$. It can prove that there is no essential difference between right- and left- preconditioning methods with respect to convergence behavior. It is worth noting that the residual vectors computed by left preconditioning and right preconditioning correspond to the preconditioned and actual residuals, respectively. The Deflated GMRES method presented reads as Algorithm \ref{CSLPDEFGMRES-1}. Furthermore, the preconditioned Generalized Conjugate Residual method (GCR) \cite{Eisenstat1983} will also be considered, as it allows for variable preconditioning inherently.

\begin{algorithm}[H]
    \caption{Preconditioned deflated GMRES for system $Au = b$}
    \label{CSLPDEFGMRES-1}
\begin{algorithmic}
	\STATE{Choose $u_0$;}
	\STATE{Compute {$r_0 =P(b-Au_0)$} and $v_1 = r_0/ \left| \left| r_0\right| \right|$;}
	\FOR{$j=1,2,...,k$ or until convergence}
		\STATE{$\tilde{v}_j=A v_j$ }
            \STATE{$w :=P \tilde{v}_j$ }
		%\STATE{\color{red}$w=M^{-1} \tilde{v}_j$ \ \% Preconditioned}
		\FOR{$i:=1,2,...,j$}
			\STATE{$h_{i,j}=\left( w, v_i\right) $}
			\STATE{$w:=w-h_{i,j}v_i$}
		\ENDFOR
		\STATE{$h_{j+1,j}:=\left| \left| w \right| \right|$}
		\STATE{$v_{j+1}=w/h_{j+1,j}$}
	\ENDFOR
	\STATE{Store $V_k=\left[ \tilde{v}_1,...,\tilde{v}_k \right] $; $H_k=\left\lbrace h_{i,j} \right\rbrace $, $1\le i \le j+1$, $1\le j \le m$}
	\STATE{Computing of minimize $y_k$ over $\left| \left| \beta e_1 - H_k y\right| \right|$ and $u_k=u_0+V_ky_k$}
	\STATE{Update approximated solution $u_k = Qb + P^Tu_k$}
 \end{algorithmic}
\end{algorithm}

\section{Parallel implementation}
\label{sec:imp}
To develop a parallel scalable iterative solver for Helmholtz problems, we intend to parallelize the deflated Krylov subspace methods in a matrix-free framework, which allows us to tackle larger-scale problems. 

The standard MPI library is employed for data communications among the processors. Based on the MPI Cartesian topology, we can partition the computational domain blockwise and allocate the variables to the corresponding processor. In domain partitioning, we carry out the partition between two grid points. Therefore, the boundary points of adjacent subdomains are adjacent grid points in the global grids. One layer of overlapping grid points is introduced outward at each interface boundary to represent the adjacent grid points. In our method, the grid unknowns are stored as an array based on the grid ordering $(i, j)$ instead of a column vector based on $x$-line lexicographic ordering.

We implement the parallel multigrid iteration based on the original global grid. According to the relationship between the fine grid and the coarse grid, the parameters of the coarse grid are determined by the grid parameters of the fine one. For example, point $(i_c,j_c)$ in the coarse grid corresponds to point $(2i_c-1, 2j_c-1)$ in the fine grid. The restriction, as well as interpolation of the grid variables, are implemented according to the stencils based on the index correspondence between the coarse and fine grid. Finally, the coarse grid problem is solved by GMRES or Bi-CGSTAB. These operations are all implemented in a matrix-free way. Currently, for a V-cycle, after reaching a manually predefined coarsest grid size, for example, each sub-domain must contain at least $2 \times 2$ grid points, the coarsening operation will stop and solve the coarse problem in parallel, which may incur some efficiency loss. 

\section{Matrix-free method} \label{sec:matrix-free}
In the Krylov subspace methods, one of the main operations is matrix-vector multiplication. As shown in Algorithm \ref{CSLPDEFGMRES-1}, we have two main matrix-vector multiplications. They are the Helmholtz operator $y=Ax$ and the deflation operator $y=Px$, where $y$ and $x$ are arbitrary variable vectors. Besides, the inverse of the CSLP preconditioner is required in ADP and TLKM.

\subsection{Helmholtz and CSLP operator}
For the Helmholtz operator, we can implement it by applying the computational stencils \cref{eq:Helmstencil} based on the finite-difference grids. For CSLP preconditioning, the preconditioner $M$ defined by (\ref{CSLP def}) has a similar computational stencil as the Helmholtz operator. Considering any grid point $(i, j)$ ($1 \le i \le nx$, $1 \le j \le ny$), define $ap$, $aw$, $ae$, $as$ and $an$ as the multipliers of $u(i, j)$, $u(i-1, j)$, $u(i+1, j)$, $u(i, j-1)$ and $u(i, j+1)$, respectively. When physical boundary conditions are encountered, it is only necessary to set the multiplier corresponding to meaningless grid points to zero. For example, if $u(i,j)$ is a left boundary grid point, $aw$ is zero. For Dirichlet boundary conditions, we can simply set $ap = 1$, $aw=ae=as=an=0$. For the Helmholtz operator, we have
\begin{equation}
    ap = \frac{4-k^2 h^2}{h^2} \quad  aw = ae = as = an = -\frac{1}{h^2} 
\end{equation} 
For the CSLP operator, we will have
\begin{equation}
    ap = \frac{4-\left(\beta_1 - \beta_2\text{i}  \right) k^2 h^2}{h^2} \quad  aw = ae = as = an = -\frac{1}{h^2} 
\end{equation}
Thus, calculating $v=A_h u$ or $y = M_h x$ can be done in a matrix-free way by Algorithm \ref{Matrix free Mat-vec}.

\begin{algorithm}[H]
    \caption{Matrix-free Matrix-Vector Multiplication $\mathbf{v}=A\mathbf{u}$.}
    \label{Matrix free Mat-vec}
    \begin{algorithmic}
    \STATE{Input: Array $\mathbf{u}$;}
    \STATE{Output: Results of $\mathbf{v} = A\mathbf{u}$;}
    \STATE{Initiate $ap$, $aw$, $ae$, $as$ and $an$;}
    \STATE{Exchange interface boundaries data;}
    
    \FOR{$j=1, 2, ..., ny$}
        \FOR{$i:=1, 2, ..., nx$}
            \IF{physical boundary grid point}
                \STATE{Reset $ap$, $aw$, $ae$, $as$ or $an$;}
            \ENDIF
            \STATE{$v(i,j)=ap*u(i,j)+ae*u(i+1,j)+aw*u(i-1,j)+an*u(i,j+1)+as*u(i,j-1)$;}
        \ENDFOR
    \ENDFOR
    \STATE{Return $v$;}
    \end{algorithmic}
\end{algorithm}

The CSLP preconditioner is approximately inverted by one parallel geometric multigrid V-cycle. Instead of the Galerkin coarsening method, we obtain the coarse-grid operator $M_{H}$ in the same way that the matrix $M_h$ is obtained on the fine mesh, which is second-order finite-difference discretization. It is the so-called re-discretization approach. The interpolation and restriction operators are implemented by applying the polynomials \cref{eq: linear iterpolation} and \cref{eq: linear restriction} based on the finite-difference grids without constructing the inter-grid transfer operators explicitly. Note that, unless mentioned, all employment of CSLP in this paper, whether on fine or coarse grids, will use the matrix-free method described in this section

\subsection{Deflation operator}
To get insights into how the deflation operator is applied on a vector, let us consider a general form that applies the entire preconditioner on a vector inside the Krylov method, that is $w=Pv$, according to the definition of the shifted deflation preconditioner, 
\begin{eqnarray}
  w &=& P_{h,\gamma}v \nonumber\\
    &=& (I_{h}-A_{h}Q_{h}+\gamma Q_{h})v \nonumber \\
    &=& v - A_{h}v' - \gamma v'  
    \label{eq: Phv-indetail}
\end{eqnarray}
where 
\begin{equation}
    v' = Q_h v = I_{2h}^{h} A_{2h}^{-1} I_{h}^{2h}v
    \label{eq: Qhv-indetail}
\end{equation}
Subsequently, we can perform the following operations to get $v'$,
\begin{subequations}
\begin{align}
\text{Restriction: } v_1 &= I_{h}^{2h} v \label{eq:matvec-restriction}\\
\text{Coarse-grid operator: } v_2 &= A_{2h}^{-1} v_1 \label{eq:coarse-grid-problem}\\
\text{Interpolation: } v'  &= I_{2h}^{h} v_2 \label{eq:matvec-interpolation}
\end{align}
\label{Qx stepbystep}
\end{subequations}
To compute \cref{eq:coarse-grid-problem}, \textit{i.e.} the coarse grid problem, we can solve $v_2$ from the following coarse system 
\begin{equation}
\label{eq: CGP}
   A_{2h} v_2 = v_1
\end{equation}
by Krylov subspace methods. Since the coarse-grid system has similar properties to the original Helmholtz system, we can also solve \cref{eq: CGP} by CSLP preconditioned methods. 

As for the TLKM preconditioner, similar to the previous exposition in (\ref{eq: Phv-indetail})-\cref{Qx stepbystep}, we need to solve $v_2$ from the following coarse system 
\begin{equation}
    \label{eq: ADEF1 CGP}
   \tilde{A}_{2h} v_2 = v_1
\end{equation}
As $\tilde{A}_{2h} = I_{h}^{2h} M_{h,(\beta_1,\beta_2)}^{-1}A_h I_{2h}^{h}$, the explicit construction of $\tilde{A}_{2h}$ requires an explicit inverse of $M_{h,(\beta_1,\beta_2)}$, which is impractical. One can also compute the main matrix-vector multiplication, e.g. $\tilde{A}_{2h} v_i = I_{h}^{2h} M_{h,(\beta_1,\beta_2)}^{-1} A_{h} I_{2h}^{h} v_i$, in the Krylov method step by step like \cref{Qx stepbystep}. Then $M_{h,(\beta_1,\beta_2)}^{-1}$ can be approximated by a multigrid cycle. However, it requires performing the fine grid Helmholtz operator $A_h$ and approximation of $M_{h,(\beta_1,\beta_2)}^{-1}$ at every Krylov iteration on the coarse level. It is computationally expensive.

To obtain a practical variant of TLKM, as suggested by \cite{erlangga2008MLKM} that $M_{h,(\beta_1,\beta_2)}^{-1} \approx (I_{h}^{2h}M_{h,(\beta_1,\beta_2)}I_{2h}^{h})^{-1}$, one can approximate $\tilde{A}_{2h}$ as follows
\begin{eqnarray}
    \tilde{A}_{2h} &=& I_{h}^{2h} M_{h,(\beta_1,\beta_2)}^{-1}A_h I_{2h}^{h} \\
    &\approx& I_{h}^{2h} I_{2h}^{h} (I_{h}^{2h}M_{h,(\beta_1,\beta_2)}I_{2h}^{h})^{-1} I_{h}^{2h} A_h I_{2h}^{h} \\
    &=& I_{h}^{2h} I_{2h}^{h} M_{2h,(\beta_1,\beta_2)}^{-1} A_{2h}
\end{eqnarray}
where $M_{2h,(\beta_1,\beta_2)}$ and $A_{2h}$ are the coarse-level operators based on the Galerkin coarsening method. Thus, we can obtain $v_2$ by solving the following coarse-grid problem
\begin{equation}
\label{eq: TLKM CGP}
   I_{h}^{2h} I_{2h}^{h} M_{2h,(\beta_1,\beta_2)}^{-1} A_{2h} v_2 = v_1
\end{equation}

One can observe that besides the inter-grid transfer operators and the CSLP preconditioner, the main matrix-vector multiplication of the coarse-grid problems \cref{eq: ADEF1 CGP} and \cref{eq: TLKM CGP} is $y=A_{2h}x$, where $x$ and $y$ denote arbitrary variable vectors on the coarse grid. Generally, the coarse grid operator $A_{2h}$ can be built by Galerkin coarsening. However, it does not fit the matrix-free framework. In the following, we will discuss several approaches for implementing this coarse grid operator.

\subsubsection{Straightforward Galerkin coarsening operator (str-Glk)} One can compute the Galerkin coarsening operator step by step, $\hat{v}_i = A_{2h} v_i$ for example,
\begin{subequations}
\begin{align}
\text{Interpolation: } v_{i1} &= I_{2h}^{h} v_i \label{eq:matvec-interpolation2}\\
\text{Fine-grid operator: } v_{i2} &= A_{h} v_{i1} \label{eq:fine-grid-problem}\\
\text{Restriction: }  \hat{v}_i  &= I_{h}^{2h} v_{i2} \label{eq:matvec-restriction2}
\end{align}
\label{A2hx stepbystep}
\end{subequations}
However, in this way, we need to perform the fine-grid Helmholtz operator $A_h$ at every Krylov iteration on the coarse level. This will be very computationally expensive. 

\subsubsection{Galerkin coarsening stencil operation (stcl-op-Glk)} According to the definition, the coarse-grid operator $A_{2h}$ is obtained by $I_h^{2h} A_h I_{2h}^h$. Since the deflation vectors and Helmholtz operator have the so-called stencils, one can perform the stencil operations to obtain the results of $I_h^{2h} A_h I_{2h}^h$ instead of constructing the matrices and carrying out the matrix-matrix multiplication. 

\begin{remark}[Stencil Operations]\label{rmk:stcl-op}
  Suppose $A \in \mathbb{C}^{n \times n}$ is symmetric and can be represented by the stencil
  \begin{equation}
      \left[\begin{array}{*{20}{c}}
            a_l&a_0&a_r\\
        \end{array}\right], \nonumber 
  \end{equation}
  $B \in \mathbb{R}^{m \times n}$ with stencil
  \begin{equation}
      \left[\begin{array}{*{20}{c}}
            b_l&b_0&b_r\\
        \end{array}\right], \nonumber 
  \end{equation}
  and $C=BA$. Then the resulting stencil of $C$ can be obtained by the following stencil operation,
  \begin{equation}
    \begin{aligned}
        &[C] = [A] \boxplus [B] \\
        &= \left[\begin{array}{*{20}{c}}
            a_l&a_0&a_r\\
            \end{array}\right] \boxplus
            \left[\begin{array}{*{20}{c}}
            b_l&b_0&b_r\\
        \end{array}\right]\\
       &=\left[\begin{array}{*{20}{c}}
            a_l b_r& a_l b_0+a_0 b_r& a_l b_l +a_0 b_0+a_r b_r &a_0 b_l +a_r b_0&a_r b_l\\
            \end{array}\right]
    \end{aligned} \nonumber
\end{equation}
\end{remark}
where $\boxplus$ is denoted as the stencil operation symbolically.

With the stencil operation, we can calculate the stencil of $A_{2h}$ as follows
\begin{equation}
\label{eq:stencil-operation-of-Helmholtz}
    \begin{aligned}
        & [A_{2h}] = [I_h^{2h}] \boxplus [A_h] \boxplus [I_{2h}^h] = \frac{1}{64} \cdot \frac{1}{h^2} \cdot \frac{1}{64} \cdot\\
        & \left[ {\begin{array}{*{20}{c}}
            1&4&6&4&1\\
            4&16&24&16&4\\
            6&24&36&24&6\\
            4&16&24&16&4\\
            1&4&6&4&1
            \end{array}} \right]_h^{2h} \hspace{-1.1em} \boxplus
            \left[ {\begin{array}{*{20}{c}}
            0&{ - 1}&0\\
            { - 1}&{4-k^2(i,j)h^2}&{ - 1}\\
            0&{ - 1}&0
            \end{array}} \right]_h \hspace{-0.8em} \boxplus
            \left[ {\begin{array}{*{20}{c}}
            1&4&6&4&1\\
            4&16&24&16&4\\
            6&24&36&24&6\\
            4&16&24&16&4\\
            1&4&6&4&1
            \end{array}} \right]_{2h}^{h}.
    \end{aligned}
\end{equation}
When encountering the boundaries, stencils such as \cref{eq:Sommerfeldstencil1} and \cref{eq:Sommerfeldstencil2} will be used. 

An insight into this method shows that it is a variant of the straightforward Galerkin coarsening approach. The convergence results of this variant are supposed to be equivalent to the straightforward Galerkin coarsening approach. Since we do not construct any matrices and store the results of $A_{2h}=I_h^{2h}A_hI_{2h}^h$, the stencil operation $[I_h^{2h}] \boxplus [A_h] \boxplus [I_{2h}^h]$ is conducted every time we need to compute $y=A_{2h}x$. One can find that the stencil operation \cref{eq:stencil-operation-of-Helmholtz} requires around $7 \times 7 \times 5 + 5 \times 5 \times 25 = 870$ multiplications. The computational cost will be prohibitive, which can be verified from the complexity analysis in the following Section \ref{sec:Complexity_analysis}.

\subsubsection{Re-discretization approach (ReD)} An alternative is to obtain the operator $A_{2h}$ by re-discretizing the Helmholtz operator on the coarse grid. Although it leads to the consequence that $P_h$ is no longer a projection, this will not break the convergence in practical application. The natural way to re-discretize is to use the same discretization method as on the fine grid, that is, use the second-order finite difference scheme to discretize the Laplace operator, and then subtract the diagonal matrix of wavenumbers. In this way, we will get the same computational stencil as \cref{eq:Helmstencil} (denoted as \textbf{ReD-$\mathcal{O}$2}). 

In addition, inspired by high-order deflation vectors \cite{dwarka2020scalable}, one can consider discretizing the Laplace operator using a fourth-order finite difference scheme (\textbf{ReD-$\mathcal{O}$4}) for the internal grid points $(x_i,y_j)$, i.e.,
\begin{equation}
    \frac{-u_{i-2,j}+16u_{i-1,j}-30u_{i,j}+16u_{i+1,j}-u_{i+2,j}}{12 h^2}=\frac{\partial^2 u}{\partial x^2}(x_i,y_j)+\mathcal{O}(h^4),
\end{equation}
or even a sixth-order finite difference scheme (\textbf{ReD-$\mathcal{O}$6}) 
\begin{equation}
    \begin{aligned}
        &\frac{2u_{i-3,j}-27u_{i-2,j}+270u_{i-1,j}-490u_{i,j}+270u_{i+1,j}-27u_{i+2,j}+2u_{i+3,j}}{180 h^2}\\
        &=\frac{\partial^2 u}{\partial x^2}(x_i,y_j) +\mathcal{O}(h^6)
    \end{aligned}
\end{equation}
where $u_{i,j}$ denotes the value of function $u(x,y)$ on the grid point $(x_i,y_j)$ (and similar for $y$-derivative). The complexity of the computational stencils obtained by using the re-discretization approach will be less than that of stencils from Galerkin coarsening. For example, if a fourth-order finite difference scheme is used, the stencil for the coarse-grid operator becomes
\begin{equation}\label{eq:Helmstencil Oh4}
\left[ {A_{2h}} \right] = \frac{1}{12 \cdot (2h)^2}\left[ {\begin{array}{*{20}{c}}
    0&0&1&0&0\\
    0&0&-16&0&0\\
    1&-16&60-{k^2_{2h}(i_c,j_c)}{(2h)^2}&-16&1\\
    0&0&-16&0&0\\
    0&0&1&0&0\\
    \end{array}} \right]_{2h}
\end{equation}
where $(i_c,j_c) \in \Omega^{2h}$.

Regarding the boundary, since the stencil is extended to $5\times5$ or even $7\times7$, more boundary grid points need to be considered. In this paper, the stencils for the boundary grid points will be changed to that of the second-order re-discretization method, for instance, the Sommerfeld boundary condition \cref{eq:Sommerfeldstencil1} and \cref{eq:Sommerfeldstencil2}. Additional experiments show that a higher-order re-discretization scheme for the boundaries will not make any improvement in the convergence compared to the second-order method.

\subsubsection{Re-discretized scheme from Galerkin coarsening (ReD-Glk)} \label{sec:ReD-Glk}
In terms of the interior grid points, the result of $I_h^{2h} A_h I_{2h}^h$ is possible to be represented by a so-called stencil notation. Combining the idea of the Galerkin coarsening and re-discretization approach, we propose to use the stencils from the results of the Galerkin coarsening operator as the finite-difference scheme to perform re-discretization on the coarse grid. 

Consider using the higher-order deflation vectors $I_h^{2h}$ and $I_{2h}^h$, for the Laplace operator, we have
\begin{equation}
    \begin{aligned}
        & [I_h^{2h}] \boxplus [-\Delta_h] \boxplus [I_{2h}^h] = \\
        & \frac{1}{64}\left[ {\begin{array}{*{20}{c}}
            1&4&6&4&1\\
            4&16&24&16&4\\
            6&24&36&24&6\\
            4&16&24&16&4\\
            1&4&6&4&1
            \end{array}} \right]_h^{2h} \hspace{-0.8em} \boxplus
            \frac{1}{{{h^2}}}\left[ {\begin{array}{*{20}{c}}
            0&{ - 1}&0\\
            { - 1}&{4}&{ - 1}\\
            0&{ - 1}&0
            \end{array}} \right]_h \hspace{-0.8em} \boxplus
            \frac{1}{64}\left[ {\begin{array}{*{20}{c}}
            1&4&6&4&1\\
            4&16&24&16&4\\
            6&24&36&24&6\\
            4&16&24&16&4\\
            1&4&6&4&1
            \end{array}} \right]_{2h}^{h} 
    \end{aligned} \nonumber
\end{equation}
\begin{equation}
\label{eq: Glk Laplace}
    \begin{aligned}
        & \Rightarrow [-\Delta_{2h}] = \frac{1}{(2h)^2}\cdot\frac{1}{256}\left[ {\begin{array}{*{20}{c}}
            -3 & -44 & -98 & -44 & -3\\
            -44& -112& 56  & -112& -44\\
            -98& 56  & 980 & 56  & -98\\
            -44& -112& 56  & -112& -44\\
            -3 & -44 & -98 & -44 & -3
            \end{array}} \right]_{2h}
    \end{aligned}
\end{equation}
where $\boxplus$ indicates the operation between stencils. Taylor series expansion analysis shows that the finite-difference scheme based on the stencil \cref{eq: Glk Laplace} is a second-order approximation of 2D Laplacian, that is,
\begin{equation}
    -\Delta_{2h}u_{2h} = 4 \cdot \{-\frac{\partial^2u}{\partial x^2}-\frac{\partial^2u}{\partial y^2}-(\frac{13}{48}\frac{\partial^4u}{\partial x^4}+\frac{1}{2}\frac{\partial^4u}{\partial x^2\partial y^2}+\frac{13}{48}\frac{\partial^4u}{\partial y^4})(2h)^2 \}+\mathcal{O}(h^4).
\end{equation}
Note that the higher-order deflation vectors enlarge the size of the stencil but maintain the accuracy of finite-difference discretization.

For the diagonal wavenumber matrix, a similar derivation can be made
\begin{equation}
    \begin{aligned}
        & [I_h^{2h}] \boxplus [ \mathcal{I}(k_{i,j}^2)_h] \boxplus [I_{2h}^h] = \\
        & \frac{1}{64}\left[ {\begin{array}{*{20}{c}}
            1&4&6&4&1\\
            4&16&24&16&4\\
            6&24&36&24&6\\
            4&16&24&16&4\\
            1&4&6&4&1
            \end{array}} \right]_h^{2h} \hspace{-0.8em} \boxplus
            \left[ {\begin{array}{*{20}{c}}
            0&0&0\\
            0&k(i,j)^2&0\\
            0&0&0
            \end{array}} \right]_h \hspace{-0.8em} \boxplus
            \frac{1}{64}\left[ {\begin{array}{*{20}{c}}
            1&4&6&4&1\\
            4&16&24&16&4\\
            6&24&36&24&6\\
            4&16&24&16&4\\
            1&4&6&4&1
            \end{array}} \right]_{2h}^{h} ,
    \end{aligned} \nonumber
\end{equation}
where $(i,j) \in \Omega^h$. Note that, for non-constant wavenumber, the result will be a $5\times5$ stencil, of which each element contains the wavenumber on up to 25 fine grid points. This is complicated to implement and leads to extra work.
% For example, the central coefficient value of the resulting stencil $[ \mathcal{I}(k_{i_c,j_c}^2)_{2h}]$ is equal to 
% \begin{equation}
%     \begin{aligned}
%         &\frac{1}{64^2} \left[ k_h(2i_c-3,2j_c+1)^2 + 16k_h(2i_c-2,2j_c+1)^2 + 36k_h(2i_c-1,2j_c+1)^2  \right.\\
%         & + 16k_h(2i_c,2j_c+1)^2 + k_h(2i_c+1,2j_c+1)^2 \\
%         & + 16k_h(2i_c-3,2j_c)^2 + 256k_h(2i_c-2,2j_c)^2 + 576k_h(2i_c-1,2j_c)^2 \\  
%         & + 256k_h(2i_c,2j_c)^2 + 16k_h(2i_c+1,2j_c)^2 \\
%         & + 36k_h(2i_c-3,2j_c-1)^2 + 576k_h(2i_c-2,2j_c-1)^2 + 1296k_h(2i_c-1,2j_c-1)^2 \\
%         & + 576k_h(2i_c,2j_c-1)^2 + 36k_h(2i_c+1,2j_c-1)^2 \\
%         & + 16k_h(2i_c-3,2j_c-2)^2 + 256k_h(2i_c-2,2j_c-2)^2 + 576k_h(2i_c-1,2j_c-2)^2 \\
%         & + 256k_h(2i_c,2j_c-2)^2 + 16k_h(2i_c+1,2j_c-2)^2 \\
%         & +  k_h(2i_c-3,2j_c-3)^2 + 16k_h(2i_c-2,2j_c-3)^2 + 36k_h(2i_c-1,2j_c-3)^2 \\
%         & + \left. 16k_h(2i_c,2j_c-3)^2 + k_h(2i_c+1,2j_c-3)^2 \right],
%     \end{aligned} \nonumber
% \end{equation}
% where $(i_c,j_c) \in \Omega^{2h}$.
In order to obtain a simple stencil, suppose the wavenumber is a constant $k$, we will have
\begin{equation}
\label{eq: ReD-Glk Kh}
    \begin{aligned}
        & [ \mathcal{I}(k^2)_{2h}] =\frac{k^2}{64^2}\left[ {\begin{array}{*{20}{c}}
            1 & 28   & 70   & 28   & 1\\
            28& 784  & 1960 & 784  & 28\\
            70& 1960 & 4900 & 1960 & 70\\
            28& 784  & 1960 & 784  & 28\\
            1 & 28   & 70   & 28   & 1\\
            \end{array}} \right]_{2h}
    \end{aligned}
\end{equation}
For non-constant wavenumber, the stencil needs to be slightly modified. For each node of the stencil, the wavenumber on the corresponding coarse grid point is used. Then the stencil of the coarse-grid operator from Galerkin coarsening is obtained by $[A_{2h}]=[-\Delta_{2h}]-[ \mathcal{I}(k_{i_c,j_c}^2)_{2h}]$.  

Since they are $5\times$5 stencils, one needs to consider two grid points on the boundary. Using the standard second-order discretization on both grid points is an easily implementable approach. An alternative method is to use the ghost grid point. Specifically, we can calculate the value of a ghost point using the method described in Equation \cref{eq:GhostPoint} for Sommerfeld boundary conditions. As for Dirichlet boundary conditions, one can approximate the value of the ghost point $u_{0,j}$ by the following relationship 
\begin{equation}
    \label{eq:DirchletGhostPoint}
    u_{1,j} =\frac{u_{0,j}+u_{2,j}}{2}
\end{equation}
Then, the aforementioned $5\times 5$ stencils will be applicable on the second grid point of the boundary. It should be noted that one can set the wavenumber of the ghost point to zero without making any approximation, which is needed by \cref{eq: ReD-Glk Kh}. The first point of the boundary can continue to be discretized using the second-order method mentioned in Section \ref{sec:Discretization}. Let us denote these two processing methods as \textbf{ReD-Glk1} and \textbf{ReD-Glk2}, respectively. However, these simplifications may result in an extra number of iterations although it may keep the wavenumber-independent convergence, which will be reflected in the numerical results in the next section.

\subsubsection{Fourth-order compact re-discretization (ReD-cmp$\mathcal{O}$4)}
Compared to the re-discretization using the explicit central finite-difference scheme, we can observe that the re-discretization scheme from Galerkin coarsening has a full $5 \times 5$ stencil for both the Laplace operator and the item related to the wavenumber. In contrast, ReD-$\mathcal{O}$2 and ReD-$\mathcal{O}$4 do not fully use the entire template and the wavenumber information of adjacent grid points. A larger stencil often means more operations and more communications in parallel computation. From ReD-Glk, one may think of a compact finite-difference approximation of the Helmholtz operator that has the following $3 \times 3$ stencil at grid point $(i_c,j_c)$:

\begin{equation}
\label{eq:GeneralFormReD}
    [A_{2h}]=\frac{1}{(2h)^2}\left[ {\begin{array}{*{20}{c}}
            a_c & a_s & a_c \\
            a_s & a_0 & a_s \\
            a_c & a_s & a_c \\
            \end{array}} \right] - \left[ {\begin{array}{*{20}{c}}
            b_c & b_s & b_c \\
            b_s & b_0 & b_s \\
            b_c & b_s & b_c \\
            \end{array}} \right] k^2
\end{equation}
where the coefficients $a$ and $b$ are non-zeros. 

The fourth-order finite difference scheme introduced by Singer and Turkel \cite{singer1998high} has a similar form as Eq. \cref{eq:GeneralFormReD}. From \cite{singer1998high}, we have
\begin{eqnarray}
    &&a_0 = \frac{10}{3},\quad a_s=-\frac{2}{3},\quad a_c=-\frac{1}{6}, \nonumber\\
    &&b_0 = \frac{2}{3}+\frac{\gamma}{36},\quad b_s=\frac{1}{12}-\frac{\gamma}{72},\quad b_c=\frac{\gamma}{144}, \nonumber
\end{eqnarray}
with $\gamma$ an arbitrary constant. Here we will use $\gamma = 1$. As for the boundary, we can use the fourth-order approximations of the Sommerfeld boundary condition given by \cite{erlangga2012iterative}. For instance, suppose $u_{0,j}$ is a ghost point on the left of $u_{1,j}$, it can be approximated by
\begin{equation}
    u_{0,j}=2 \text{i} k h(1-\frac{k^2h^2}{6}) u_{1,j}+u_{2,j}.
\end{equation}
For the ghost corner point $u_{0,0}$, we approximate it by $u_{0,0}=\frac{u_{1,0}+u_{0,1}}{2}$.

A potential way to construct a nine-point compact re-discretization scheme by aligning the smallest eigenvalues refers to Appendix \ref{Appendix B}.

%In this paper, the predefined coarsest global grid size is $ nc_x \times nc_y = 9 \times 9$ as the maximum number of processors we use is $4 \times 4$. 

% Our analysis leads to the algorithm in \cref{alg:buildtree}.

% \begin{algorithm}
% \caption{Build tree}
% \label{alg:buildtree}
% \begin{algorithmic}
% \STATE{Define $P:=T:=\{ \{1\},\ldots,\{d\}$\}}
% \WHILE{$\#P > 1$}
% \STATE{Choose $C^\prime\in\mathcal{C}_p(P)$ with $C^\prime := \operatorname{argmin}_{C\in\mathcal{C}_p(P)} \varrho(C)$}
% \STATE{Find an optimal partition tree $T_{C^\prime}$ }
% \STATE{Update $P := (P{\setminus} C^\prime) \cup \{ \bigcup_{t\in C^\prime} t \}$}
% \STATE{Update $T := T \cup \{ \bigcup_{t\in\tau} t : \tau\in T_{C^\prime}{\setminus} \mathcal{L}(T_{C^\prime})\}$}
% \ENDWHILE
% \RETURN $T$
% \end{algorithmic}
% \end{algorithm}

\section{Experimental results}
\label{sec:experiments}
The numerical experiments are conducted on the Linux supercomputer DelftBlue \cite{DHPC2022}. DelftBlue runs on the Red Hat Enterprise Linux 8 operating system. Each compute node is equipped with two Intel Xeon E5-6248R processors with 24 cores at 3.0 GHz, 192 GB of RAM, a memory bandwidth of 132 GByte/s per socket, and a 100 Gbit/s InfiniBand card. In our experiments, the solver is developed in Fortran 90. On DelftBlue, the code is compiled using GNU Fortran 8.5.0 with the compiler options \verb|-O3| for optimization purposes. Open MPI library (version 4.1.1) is employed for message passing.

In the numerical experiments, the convergence test is based on the $l_2$-norm of the residual. Unless mentioned, the preconditioned relative residual is reduced to $10^{-6}$ by the deflated GMRES algorithm \ref{CSLPDEFGMRES-1}. Since it is a two-level method, it may still be  expensive to solve the coarse-grid problem directly. The (CSLP-preconditioned) GMRES is employed for approximating the inverse of the coarse-grid operator. The (preconditioned) relative residuals on the coarse level are reduced to a certain tolerance. In our study of the convergence properties of the deflation-preconditioned Krylov-subspace method for solving the Helmholtz problems, we aimed to obtain the approximation of $A_{2h}^{-1}$ as accurately as possible. To achieve this, we set the tolerance for the coarse-grid solver to $10^{-12}$. Furthermore, we will explore the optimal tolerance for the coarse-grid solver in our study.

The inverse of the CSLP preconditioner is approximated by a multigrid V-cycle, where full-GMRES is used to reduce the relative residuals to $10^{-8}$ on the predefined coarsest level \cite{jchen2D2022}.

%or the number of matrix-vector multiplications (denote as \#Matvec)
This section will mainly illustrate the numerical performance of our solver for the model problems on DelftBlue. The number of iterations (denote as \#iter) will be the main metric to estimate the convergence. The speedup and parallel efficiency are used to study the parallel performance. The Wall-clock time for the preconditioned Krylov solver to reach the stopping criterion is denoted as $t_w$. The speedup $S_p$ is defined by
\begin{equation}
    S_p=\frac{t_{w, r}}{t_{w, p}}
\end{equation}
where $t_{w, r}$ is the Wall-clock time for the reference case and $t_{w, p}$ is the Wall-clock time for parallel computation. The parallel efficiency $E_p$ is given by 
\begin{equation}
    E_P=\frac{S_p}{np/{np}_r}=\frac{t_{w, r} \cdot {np}_r }{t_{w, p} \cdot np}
\end{equation}
where $np$ ($np_r$) is the (reference) number of processors.

\subsection{Complexity analysis of the coarse grid operators}\label{sec:Complexity_analysis}
Recalling from the previous section, there are the following possible coarse-grid operators in matrix-free:
\begin{itemize}
    \item \textbf{str-Glk}: Straightforward Galerkin coarsening approach
    \item \textbf{stcl-op-Glk}: Galerkin coarsening stencil operation
    \item \textbf{ReD-$\mathcal{O}$2}: Re-discretization using the second-order finite-difference scheme
    \item \textbf{ReD-$\mathcal{O}$4}: Re-discretization using the fourth-order finite-difference scheme
    \item \textbf{ReD-$\mathcal{O}$6}: Re-discretization using the sixth-order finite-difference scheme
    \item \textbf{ReD-cmp$\mathcal{O}$4}: Re-discretization using a scheme derived from the fourth-order compact finite difference of the Helmholtz equation
    \item \textbf{ReD-Glk1}: Re-discretization using a scheme derived from the Galerkin coarsening approach (ReD-$\mathcal{O}$2 are used for two boundary grid points)
    \item \textbf{ReD-Glk2}: Re-discretization using a scheme derived from the Galerkin coarsening approach (A ghost point is included for the second boundary grid point and ReD-$\mathcal{O}$2 is only for the first boundary grid point)
\end{itemize}

To give a preliminary comparison of the above methods for obtaining coarse grid operators, we conduct a rough complexity analysis to quantify the FLOPs for performing the coarse grid operator $y=A_{2h}x$ by different methods. 

\begin{remark}[FLOPs of sparse matrix-vector multiplication]\label{rmk:Flop-matvec}
    Given $x \in \mathbb{R}^{p}$ and sparse $A \in \mathbb{R}^{p\times p}$, the number of non-zero elements ($nnz$) of each row is $q$, then the upper bound of total FLOPs for $Ax$ is $2pq$
\end{remark}

Suppose the number of the fine grid points is $N$ and that of the coarse grid is $M$, if the matrices are constructed explicitly, we will have $A_h \in \mathbb{R}^{N \times N}$, $A_{2h} \in \mathbb{R}^{M \times M}$, $I_h^{2h} \in \mathbb{R}^{M \times N}$ and $I_{2h}^{h} \in \mathbb{R}^{N \times M}$. In Table \ref{tab:complexity-CGP}, we list the $nnz$ of the relevant matrices and the approximate FLOPs of $A_{2h}x$ for different methods. Since $4M \approx N$, the total FLOPs of the straightforward Galerkin coarsening operator is around $162M$. Note that the extra FLOPs required to perform the stencil operation \cref{eq:stencil-operation-of-Helmholtz} will be up to $870\times2=1740$. This is very expensive. Since the boundaries are not considered,  ReD-Glk represents both ReD-Glk1 and ReD-Glk2 in Table \ref{tab:complexity-CGP}. To verify the results of the complexity analysis, we also test the elapsed CPU time for performing the coarse grid operators derived from Galerkin coarsening approach. One can find the ratios between the elapsed CPU time are fairly consistent with the estimated FLOPs. From the complexity analysis of the coarse grid operator, without considering the overall convergence characteristics, the re-discretization approach seems to be preferable in the frame of matrix-free implementation.
\begin{table}[htbp]
    \centering
    \caption{Upper bound to the number of non-zero elements of each relevant matrix and the approximate FLOPs as well as the elapsed CPU time for performing the coarse grid operator once. The coarse grid size is $369 \times 121$.}
    \scalebox{0.65}{
    \begin{tabular}{cccccccccc}
    \hline
        Methods & \multicolumn{3}{c}{str-Glk} & stcl-op-Glk & ReD-$\mathcal{O}$2 & ReD-$\mathcal{O}$4 & ReD-$\mathcal{O}$6 & ReD-cmp$\mathcal{O}$4 &ReD-Glk \\ \hline
        Operators & $I_{h}^{2h}$ &$A_h$ & $I_{2h}^h$  & $\boxplus$ &$A_{2h}$ &$A_{2h}$ &$A_{2h}$ &$A_{2h}$ &$A_{2h}$ \\
        Max. nnz per row & 25 & 5 & 9 & - & 5 & 9 & 13 & 9 & 25 \\
        FLOPs & \multicolumn{3}{c}{$50M+10N+18N$} & $(50+1740)M$ & $10M$ & $18M$ & $26M$ & $18M$ & $50M$ \\ 
        CPU time &\multicolumn{3}{c}{$4.82\times 10^{-4}$} & $7.16\times10^{-3}$ & $3.38 \times 10^{-5}$ & - & - & - & $1.80 \times 10^{-4}$ \\\hline
    \end{tabular}}
    \label{tab:complexity-CGP}
\end{table}

\subsection{Scalable convergence}
As we introduced, in the context of solving the Helmholtz problem using iterative solvers, the number of iterations required increases significantly as the wavenumber increases. A recent study by Dwarka and Vuik \cite{dwarka2020scalable} has shown promising results in achieving convergence properties that are independent of the wavenumber using high-order deflation preconditioning methods. In their work, the Galerkin coarsening approach is employed to derive the coarse-grid operator. In this subsection, we further explore various scenarios utilizing the re-discretization method for the coarse-grid operator to investigate which approach can achieve wavenumber-independent convergence. By examining different cases and comparing the outcomes, we aim to provide insights into the effectiveness of these methods.

\subsubsection{Low-order deflation} 
To observe the possible impact on the convergence behavior of using the re-discretization approach, we first compare the convergence behavior of low-order deflation using the straightforward Galerkin coarsening and second-order re-discretization approaches. As mentioned, regardless of the computational cost, we can obtain the coarse-grid operator as \cref{A2hx stepbystep} in the frame of the Galerkin coarsening approach. \Cref{tab: MP2a Galerkin and reO2 A2h no CSLP and CSLP CGP} gives the number of outer iterations required to solve MP-1a using A-DEF1 preconditioned GMRES for different $kh$. In parentheses is the approximate number of iterations required to solve the coarse-grid problem (CGP) once by (preconditioned) GMRES. One can find the number of outer iterations is consistent with the results in \cite{sheikh2013convergence} when the straightforward Galerkin coarsening approach is used. Since the coarse-grid problem has similar properties to the original Helmholtz problem, one likes to solve the coarse-grid problem combined with the CSLP preconditioner. \Cref{tab: MP2a Galerkin and reO2 A2h no CSLP and CSLP CGP} shows that the CSLP preconditioner for the coarse-grid problem solver can accelerate the convergence on the coarse level, while the outer number of iterations remains the same.
\begin{table}
    \centering
    \caption{The number of iterations required to solve MP-2a using A-DEF1/TLKM preconditioned GMRES for different $k$ and $kh$. The coarse-grid problem (CGP) obtained by str-Glk or ReD-$\mathcal{O}$2 approach is solved by full-GMRES and CSLP preconditioned GMRES (CSLP-GMRES). In parentheses is the number of iterations required to solve the coarse grid problem once.}
    \label{tab: MP2a Galerkin and reO2 A2h no CSLP and CSLP CGP}
    \scalebox{0.75}{
    \begin{threeparttable}[htbp]
        \begin{tabular}{llllllll}
        \hline
        & & & \multicolumn{4}{c}{A-DEF1} & {TLKM}   \\ 
        & & & \multicolumn{2}{c}{str-Glk} & \multicolumn{2}{c}{ReD-$\mathcal{O}$2} & {ReD-$\mathcal{O}$2}   \\ 
        Grid size & $k$ & $kh$&  (GMRES) &  (CSLP-GMRES)   &  (GMRES)      &  (CSLP-GMRES)  &  (GMRES)      \\ 
        \hline 
        33 $\times$ 33     & 20    & 0.625   & 8(55)      & 8(54)   & 21(82)       & 21(63)   & 9(56) \\
        65 $\times$ 65     & 40    & 0.625   & 16(259)    & 16(243) & 40(326)      & 41(248)  & 20(241) \\
        129 $\times$ 129   & 80    & 0.625   & 42(1069)   & 44(911) & 97($>$3000)\tnote{1}  & 95(1099) & 70(975) \\
        &&&&&&& \\
        65 $\times$ 65     & 20    & 0.3125  & 6(131)     & 6(74)   & 17(123)      & 17(59)   & 6(68) \\
        129 $\times$ 129   & 40    & 0.3125  & 8(533)     & 9(250)  & 20(660)      & 20(206)  & 9(218) \\
        257 $\times$ 257   & 80    & 0.3125  & 19($>$1500)  & 18(944) &     - \tnote{2}        &    -     & 16(827) \\
        %321 $\times$ 321   & 100   & 0.3125  &            &         &              &          & 24(1369)\\
        \hline
        \end{tabular}
        \begin{tablenotes}
            \item [1] ``$>$" indicates it does not converge to the specified residual tolerance within a certain number of iterations.
            \item [2] ``-" indicates that the numerical experiments were not conducted; this is not due to divergence, as the existing data sufficiently illustrate the behavior.
        \end{tablenotes}
    \end{threeparttable}}
\end{table}

As for the second-order re-discretization approach, \Cref{tab: MP2a Galerkin and reO2 A2h no CSLP and CSLP CGP} shows that the second-order re-discretized $A_{2h}$ leads to an increase in the outer iterations. One can observe that the behavior of the outer iterations changing with the wavenumber and $kh$ is consistent with the straightforward Galerkin coarsening approach. For example, when the wavenumber is constant, a smaller $kh$ results in smaller outer iterations.

To observe whether the relative convergence behavior between different preconditioners changes when the re-discretization approach is used in different preconditioners, one can also utilize ReD-$\mathcal{O}$2 in the coarse grid problem of practical TLKM. In \cref{tab: MP2a Galerkin and reO2 A2h no CSLP and CSLP CGP}, the practical TLKM preconditioned GMRES is used to solve MP-1a. Note that the number of iterations required to solve the coarse grid problem is a bit more than that using the CSLP preconditioned methods. This is because the coefficient matrix in the coarse-grid system in TLKM is $I_h^{2h} I_{2h}^h M_{2h, (\beta_1, \beta_2)}^{-1} A_{2h}$, while that of A-DEF1 is $M_{2h, (\beta_1, \beta_2)}^{-1} A_{2h}$ when we solve the coarse-grid problem by CSLP preconditioned Krylov methods. The number of outer iterations appears that TLKM implemented in a re-discretization way performs similarly to A-DEF1 using the straightforward Galerkin Coarsening approach. Even the second-order re-discretization approach is used for the coarse grid problem, TLKM is still a computationally expensive method since an application of $M^{-1}_{h,(\beta_1, \beta_2)}$ on the fine grid for every coarse grid iteration.

For MP-2b, we also observe the convergence behavior mentioned above. The only difference is that the introduction of the Sommerfeld boundary conditions makes it a bit easier to solve the coarse grid problem, see \cref{tab: MP2b Galerkin and reO2 A2h CSLP CGP}. 
\begin{table}[htbp]
\centering
\caption{The number of iterations required to solve MP-2b using A-DEF1/TLKM preconditioned GMRES for different $k$ and $kh$. The coarse-grid problem (CGP) obtained by str-Glk or ReD-$\mathcal{O}$2 approach is solved by CSLP preconditioned GMRES. In parentheses is the number of iterations required to solve the coarse grid problem once.}
\label{tab: MP2b Galerkin and reO2 A2h CSLP CGP}
\scalebox{0.75}{
\begin{tabular}{llllll}
\hline
& & & \multicolumn{2}{c}{A-DEF1} & {TLKM}   \\ 
& & & {str-Glk} & { ReD-$\mathcal{O}$2} & {ReD-$\mathcal{O}$2}   \\ 
Grid size & $k$ & $kh$&  (CSLP-GMRES)   &  (CSLP-GMRES)  &  (GMRES)      \\                   
\hline
33 $\times$ 33        & 20   & 0.625    & 9(42)   & 20(39)  & 9(44) \\
65 $\times$ 65        & 40   & 0.625    & 13(297) & 26(123) & 13(94) \\
129 $\times$ 129      & 80   & 0.625    & 22(258) & 41(330) & 24(250) \\
&&&&\\
65 $\times$ 65        & 20   & 0.3125   & 6(41)   & 17(35)  & 6(56) \\
129 $\times$ 129      & 40   & 0.3125   & 7(81)   & 19(76)  & 7(104) \\
257 $\times$ 257      & 80   & 0.3125   & 10(221) & 20(204) & 8(233) \\
321 $\times$ 321      & 100  & 0.3125   & 11(303) & -       & 9(306) \\
\hline
\end{tabular}}
\end{table}

As the non-constant wavenumber problem, the so-called 2D Wedge problem is considered. \Cref{tab: MP3 Galerkin and reO2 A2h CSLP CGP} confirms that our matrix-free two-level deflation methods also work for the case with a non-constant wavenumber and exhibit similar convergence behavior to the constant-wave-number case discussed above. The largest difference is that, in \cref{tab: MP3 Galerkin and reO2 A2h CSLP CGP}, we find that when the frequency is constant, a smaller $kh$ does not lead to a smaller number of outer iterations. This may also be a disadvantage of using the re-discretization approach. 

\begin{table}[htbp]
\centering
\caption{The number of iterations required to solve the 2D Wedge problem using A-DEF1/TLKM preconditioned GMRES for different frequencies $f$ and grid size. The coarse-grid problem (CGP) obtained by str-Glk or ReD-$\mathcal{O}$2 approach is solved by CSLP preconditioned GMRES. In parentheses is the number of iterations required to solve the coarse grid problem once.}
\label{tab: MP3 Galerkin and reO2 A2h CSLP CGP}
\scalebox{0.75}{
\begin{tabular}{llllll}
\hline
& & & \multicolumn{2}{c}{A-DEF1} & {TLKM}   \\ 
& & & {str-Glk} & {ReD-$\mathcal{O}$2} & {ReD-$\mathcal{O}$2}   \\ 
Grid size & $k$ & $kh$&  (CSLP-GMRES)   &  (CSLP-GMRES)  &  (GMRES)      \\             
\hline
73 $\times$ 121      & 10      & 0.34907     & 8(112)   & 24(104) & 8(167)\\
145 $\times$ 241     & 20      & 0.34907     & 9(278)   & 31(242) & 9(362)\\
289 $\times$ 481     & 20      & 0.17453     & 7(257)   & 31(245) & 7(570)\\
289 $\times$ 481     & 40      & 0.34907     & 12(629)  & 34(547) & 12(751)\\
361 $\times$ 601     & 50      & 0.34907     & 14(815)  & 29(696) & 13(948)\\
\hline
\end{tabular}}
\end{table}

\subsubsection{Higher-order deflation}
To get scalable convergence with respect to the wavenumber, this section considers using the higher-order deflation vectors, \textit{i.e.} APD.

\cref{tab: MP2b ZO4 Galerkin reO2 reO4 reO6 ReD-Glk A2h CSLP CGP} shows the number of iterations required to solve the constant wavenumber problem (MP-2b) using APD-preconditioned GMRES. We found that when the str-Glk approach is used for the coarse grid problem, the number of iterations required shows the wavenumber independence for the case with the wavenumber up to $1000$. The slight increase is due to the fact that the coarse mesh problem does not converge to sufficient tolerance within a certain number of iterations. However, when the coarse grid operator is obtained by the ReD-$\mathcal{O}$2, it requires more outer iterations, which increase significantly with the increase of the wavenumber. When using the ReD-$\mathcal{O}$4 to obtain the coarse grid operator, the number of outer iterations is a bit less than that of ReD-$\mathcal{O}$2. However, it is still wavenumber dependent. Besides, compared to ReD-$\mathcal{O}$4, a higher-order method (ReD-O6) and a fourth-order compact scheme (ReD-cmp$\mathcal{O}$4) do not yield gains anymore. But one should note that a fourth-order compact scheme for coarse grid re-discretization significantly reduces the number of iterations required to solve the coarse grid problem. As for using ReD-Glk, we can also obtain considerable wavenumber independence. But the number of outer iterations for using ReD-Glk1 is nearly double that of using str-Glk. This should be due to the simplified boundary treatment mentioned in Section \ref{sec:ReD-Glk}. In contrast, ReD-Glk2 proves to be a more favorable choice, which confirms that treatment of the boundary correlates to the number of outer iterations.

\begin{table}
\centering
\caption{The number of iterations required to solve MP-2b using APD-preconditioned GMRES for different $k$ and $kh$. The coarse-grid operator is obtained by str-Glk, ReD-$\mathcal{O}$2, ReD-$\mathcal{O}$4, ReD-$\mathcal{O}$6 and ReD-Glk, respectively. The coarse-grid problem is solved by CSLP preconditioned GMRES. In parentheses is the number of iterations required to solve the coarse grid problem once.}
\label{tab: MP2b ZO4 Galerkin reO2 reO4 reO6 ReD-Glk A2h CSLP CGP}
\scalebox{0.7}{
\begin{threeparttable}[htbp]
\begin{tabular}{llllllllll}
\hline 
Grid size & $k$ & $kh$ & str-Glk & ReD-$\mathcal{O}$2 & ReD-$\mathcal{O}$4 & ReD-cmp$\mathcal{O}$4 & ReD-$\mathcal{O}$6 & ReD-Glk1 & ReD-Glk2  \\
\hline
65 $\times$ 65    & 40  & 0.625  & 7(126)      & 20(98)      & 17(106)  & 19(91)   & -        & 12(130)  & 9(128) \\
129 $\times$ 129  & 80  & 0.625  & 7(236)      & 30(305)     & 18(298)  & 20(185)  & 20(318)  & 12(249)  & 9(251)\\
257 $\times$ 257  & 160 & 0.625  & 7(495)      & 87(731)  & 19(650)  & 23(362)  & 25(700)  & 12(541)  & 9(585)\\
513 $\times$ 513  & 320 & 0.625  & 7($>$1000)  & $>$100\tnote{1}  & 23(1330) & 28(690)  & - \tnote{2}       & 12(1068) & 10(1276)\\
%1025 $\times$ 1025& 640 & 0.625  & 8($>$2000)  & -           & -        & -        & -        & -        & -\\
%1601 $\times$ 1601& 1000& 0.625  & 9($>$2000)  & -           & -        & -        & -        & -        & -\\
&&&&&&&&\\
129 $\times$ 129  & 40  & 0.3125 & 5(142)      & 18(76)      & 18(81)   & 18(76)   & 18(85)   & 9(160)   & 7(144)\\
257 $\times$ 257  & 80  & 0.3125 & 5(228)      & 19(205)     & 18(212)  & 18(188)  & 18(219)  & 9(264)   & 7(231)\\
513 $\times$ 513  & 160 & 0.3125 & 5(411)      & 21(438)     & 18(458)  & 19(395)  & -        & 9(460)   & 7(432)\\
1025 $\times$ 1025& 320 & 0.3125 & 5(809)           & 28(885)     & 20(909)  & 20(775)  & -        & 9(904)   & 6(859) \\
2049 $\times$ 2049& 640 & 0.3125 & 5(1630)           & 53(1722)    & 23(1763) & -        & -        & -        & 6(1690) \\
%3201 $\times$ 3201& 1000& 0.3125 & -           & 99($>$2000) & -        & -        & -        & -        & \\
\hline
\end{tabular}
\begin{tablenotes}
    \item [1] ``$>$" indicates it does not converge to the specified residual tolerance within a certain number of iterations.
    \item [2] ``-" indicates that the numerical experiments were not conducted; this is not due to divergence, as the existing data sufficiently illustrate the behavior.
\end{tablenotes}
\end{threeparttable}}
\end{table}

For the non-constant wavenumber problem, the results for solving the so-called Wedge problem and  the heterogeneous Marmousi problem are given in \cref{tab: MP3 ZO4 Galerkin reO2 reO4 reO6 ReD-Glk A2h CSLP CGP} and \cref{tab: MP4 ZO4 Galerkin reO2 reO4 reO6 ReD-Glk A2h CSLP CGP}, respectively. The convergence behaviors exhibited by using different coarse grid operators are consistent with solving the constant wavenumber problem. It is worth mentioning that the stencil \cref{eq: ReD-Glk Kh} in ReD-Glk is obtained based on the constant wavenumber assumption. But it also shows wavenumber independence in the non-constant wavenumber problem. This is consistent with the work of Dwarka and Vuik \cite{dwarka2020scalable}.

\begin{table}[htbp]
\centering
\caption{The number of iterations required to solve the 2D Wedge problem using APD-preconditioned GMRES for different frequencies $f$ and grid size. The coarse-grid operator is obtained by str-Glk, ReD-$\mathcal{O}$2, ReD-$\mathcal{O}$4, ReD-$\mathcal{O}$6, ReD-cmp$\mathcal{O}$4 and ReD-Glk, respectively. The coarse-grid problem is solved by CSLP preconditioned GMRES. In parentheses is the number of iterations required to solve the coarse grid problem once.}
\label{tab: MP3 ZO4 Galerkin reO2 reO4 reO6 ReD-Glk A2h CSLP CGP}
\scalebox{0.7}{
\begin{tabular}{llllllllll}
\hline 
Grid size & $f$ & $kh$ & str-Glk & ReD-$\mathcal{O}$2 & ReD-$\mathcal{O}$4 & ReD-$\mathcal{O}$6  & ReD-cmp$\mathcal{O}$4 & ReD-Glk1 & ReD-Glk2\\
\hline
73 $\times$ 121       & 10       & 0.349     & 7(145)    & 22(104)   & 22(108)    & 22(112)   & 22(99)    & 12(155)     &9(138)\\
145 $\times$ 241      & 20       & 0.349     & 6(301)    & 28(244)   & 27(243)    & 27(248)   & 28(228)   & 12(346)     &9(303)\\
289 $\times$ 481      & 20       & 0.174     & 6(290)    & 31(245)   & 31(246)    & -         & 32(237)   & -           &-\\
289 $\times$ 481      & 40       & 0.349     & 6(580)    & 31(535)   & 29(519)    & 29(533)   & 30(491)   & 12(718)     &9(585)\\
577 $\times$ 961      & 80       & 0.349     & 6(1206)   & 37(1200)  & 30(1175)   & -         & 31(1069)  & 12(1500)    &9(1255)\\
1153 $\times$ 1921    & 160      & 0.349     & 6($>$1500)& -         & 34($>$2000)& -         & 35(2353)  & 12($>$2500) &8($>$2500) \\
\hline
\end{tabular}}
\end{table}

\begin{table}[htbp]
\centering
\caption{The number of iterations required to solve the 2D Marmousi problem using APD-preconditioned GMRES for different frequencies $f$ and grid size. The coarse-grid operator is obtained by str-Glk, ReD-$\mathcal{O}$2, ReD-$\mathcal{O}$4 and ReD-Glk, respectively. The coarse-grid problem is solved by CSLP preconditioned GMRES. In parentheses is the number of iterations required to solve the coarse grid problem once.}
\label{tab: MP4 ZO4 Galerkin reO2 reO4 reO6 ReD-Glk A2h CSLP CGP}
\scalebox{0.75}{
\begin{tabular}{llllllll}
\hline 
Grid size & $f$ & $kh$ & str-Glk & ReD-$\mathcal{O}$2 & ReD-$\mathcal{O}$4  & ReD-Glk1 & ReD-Glk2\\
\hline
737 $\times$ 241    & 10   & 0.5236   & 7(775)      & 38(748)   & 30(762)     &13(851)  &10(802) \\
1473 $\times$ 481    & 20   & 0.5236   & 7(1858)     & 71(1988)  & 34(1947)    &13(1995) &10(1923)  \\
2945 $\times$ 961    & 40   & 0.5236   & 7($>$2500)  & -         & 50($>$2500) &13($>$2500) &11($>$2500)\\
\hline
\end{tabular}}
\end{table}

The complexity analysis in Section \ref{sec:Complexity_analysis} shows that in terms of performing coarse grid operator one time, the re-discretization approach is preferable. However, comparing the Wall-clock time required to solve the Marmousi problem, see \cref{tab: cpuTIME MP4 ZO4 Galerkin reO2 reO4 reO6 ReD-Glk A2h CSLP CGP}, we found that str-Glk is the best choice in terms of time consumption. The reason is, for the two-level method in 2D, the benefit of str-Glk in the number of outer iterations is much greater than the extra cost of the coarse grid operator. But the wavenumber independence exhibited by ReD-Glk suggests that it might be an applicable choice if our method is extended to multi-level.

\begin{table}[htbp]
\centering
\caption{The Wall-clock time required to solve the 2D Marmousi problem using APD-preconditioned GMRES for different frequencies $f$ and grid size with 12 processors. The coarse-grid operator is obtained by str-Glk, ReD-$\mathcal{O}$2, ReD-$\mathcal{O}$4 and ReD-Glk, respectively. The coarse-grid problem is solved by CSLP preconditioned GMRES.}
\label{tab: cpuTIME MP4 ZO4 Galerkin reO2 reO4 reO6 ReD-Glk A2h CSLP CGP}
\scalebox{0.75}{
\begin{tabular}{lllllllll}
\hline 
Grid size & $f$ & np & str-Glk & ReD-$\mathcal{O}$2 & ReD-$\mathcal{O}$4  & ReD-Glk1 & ReD-Glk2\\
\hline
737 $\times$ 241      & 10       & 12     & 120.67       & 615.54   & 402.79       & 227.18  &164.79\\
1473 $\times$ 481      & 20       & 12     & 1039.96      & 9606.14  & 4386.57      & 1834.36  &1376.40  \\
%2945$\times$ 961     & 40       & 48     & 3122.36      & -        & 17700.24     & 528742   & 4388.87\\
\hline
\end{tabular}}
\end{table}

\subsection{Tolerance for coarse-grid problem}
It appears that using the ReD-Glk2 format for coarse-grid re-discretization leads to near wavenumber-independent convergence. Despite this, a large number of iterations is still needed for the coarse-grid problem, especially for large wavenumbers. This is partially due to the fact that we are currently employing a two-level deflation method. Furthermore, as mentioned earlier, a strict tolerance is used for the convergence test of the coarse-grid solver. Consequently, this section explores the range of possible values for the tolerance of the coarse-grid solver. 

Using higher-order deflation vectors and ReD-Glk2 for re-discretization of the coarse-grid system, the present APD-preconditioned GMRES is utilized to solve the linear system ($A_h u_h = b_h$), reducing the preconditioned relative residual to $10^{-6}$. The stopping criteria for the coarse-grid iterative solver ($A_{2h} v_2 = v_1$) ranged from $10^{-2}$ to $10^{-13}$. Table \ref{tab:MP2BC2_outer_iter_varies_CoarseTol} shows that, for various model problems, the number of outer iterations required remains constant for all convergence criteria for the coarse grid.

\begin{table}[htbp]
\centering
\caption{The number of APD-preconditioned GMRES iterations required to solve corresponding model problems when using different stopping criteria for the coarse-grid iterative solver}
\scalebox{0.65}{
\begin{tabular}{ccccccccccccc}
\hline
Model Problems & $10^{-2}$ & $10^{-3}$ & $10^{-4}$ & $10^{-5}$ & $10^{-6}$ & $10^{-7}$ & $10^{-8}$ & $10^{-9}$ & $10^{-10}$ & $10^{-11}$ & $10^{-12}$ & $10^{-13}$ \\ \hline
MP-2b, $k=320$, $513 \times 513$   & 10 & 10 & 10 & 10 & 10 & 10 & 10 & 10 & 10 & 10 & 10 & 10\\
MP-2b, $k=320$, $1025 \times 1025$ & 7  & 7  & 6  & 6  & 6  & 6  & 6  & 6  & 6  & 6  & 6  & 6 \\
Wedge, $f=\SI{40}{\hertz}$, $289 \times 481$ & 9  & 9  & 9  & 9  & 9  & 9  & 9  & 9  & 9  & 9 & 9 & 9  \\
Marmousi, $f=\SI{20}{\hertz}$, $1473 \times 481$  & 11 & 10 & 10 & 10 & 10 & 10 & 10 & 10 & 10 & 10 & 10 & 10    \\ \hline
\end{tabular}}
\label{tab:MP2BC2_outer_iter_varies_CoarseTol}
\end{table}

\begin{figure}
    \centering
    \includegraphics[width=0.45\textwidth]{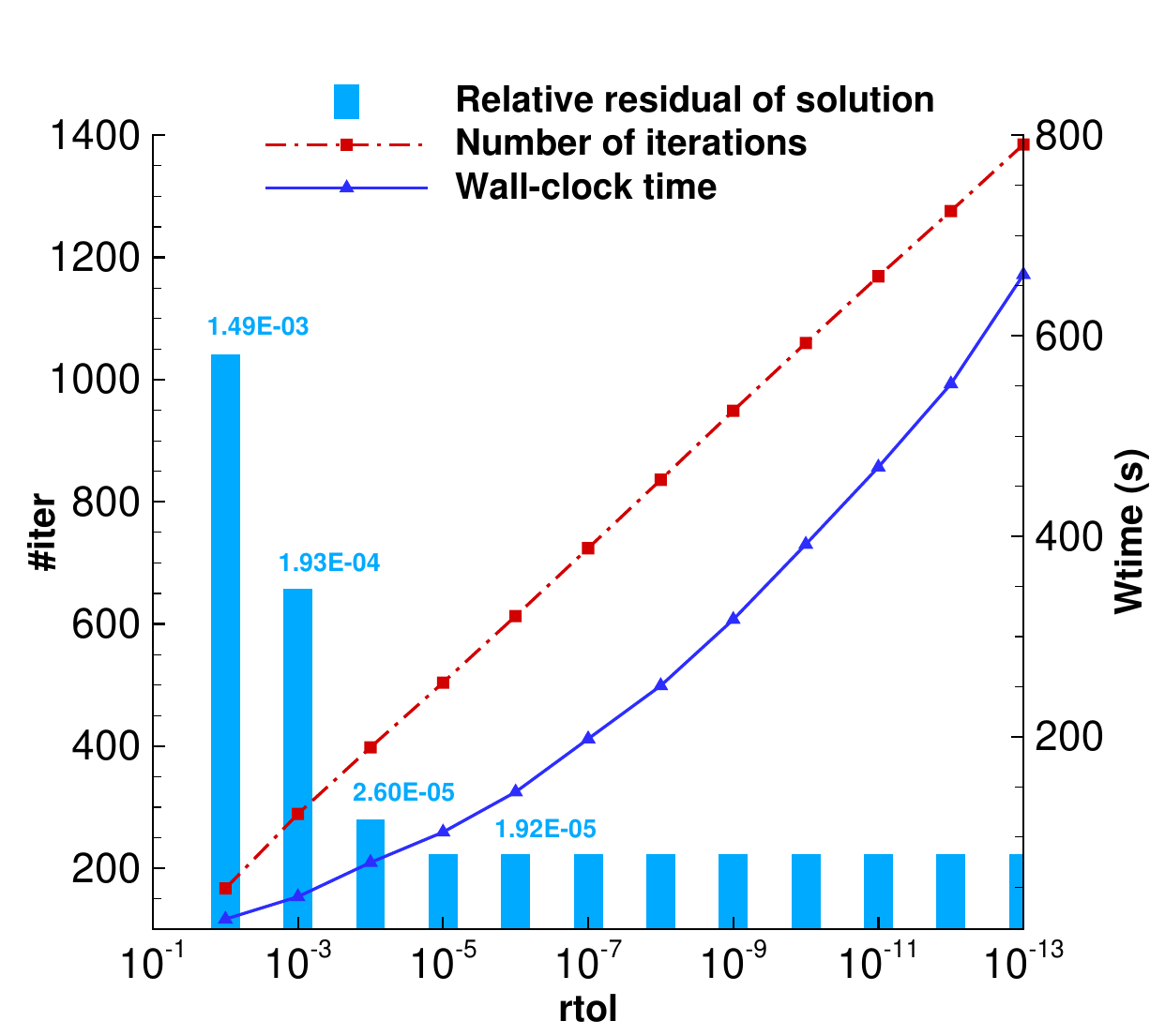}
    \caption{Tolerance of the coarse-grid solver for MP-2b with wavenumber $k=320$ and grid size $513 \times 513$. The present APD-preconditioned GMRES is employed for outer iteration.}
    \label{fig:MP2BC2_on_CoarseSolve_rtol}
\end{figure}

For the constant wavenumber problem, MP-2b with wavenumber $k=320$ and grid size $513 \times 513$, Figure \ref{fig:MP2BC2_on_CoarseSolve_rtol} displays the change in iterations required to solve a single coarse-grid problem with the CSLP-preconditioned GMRES solver, for various tolerance, along with the total Wall-clock time and the relative residual of the final solution. Note that the number of iterations for the CSLP-preconditioned GMRES solver to solve a single coarse-grid problem increases linearly as the size of the tolerance decreases. The Wall-clock time elapsed increases even faster as the tolerance becomes more stringent. Once the tolerance exceeds $10^{-5}$, the relative residual of the final solution increases, while it remains constant for the tolerance less than $10^{-6}$. This pattern holds true for larger grid sizes for MP-2b. 

\begin{figure}
     \centering
     \begin{subfigure}[b]{0.45\textwidth}
         \centering
         \includegraphics[width=\textwidth]{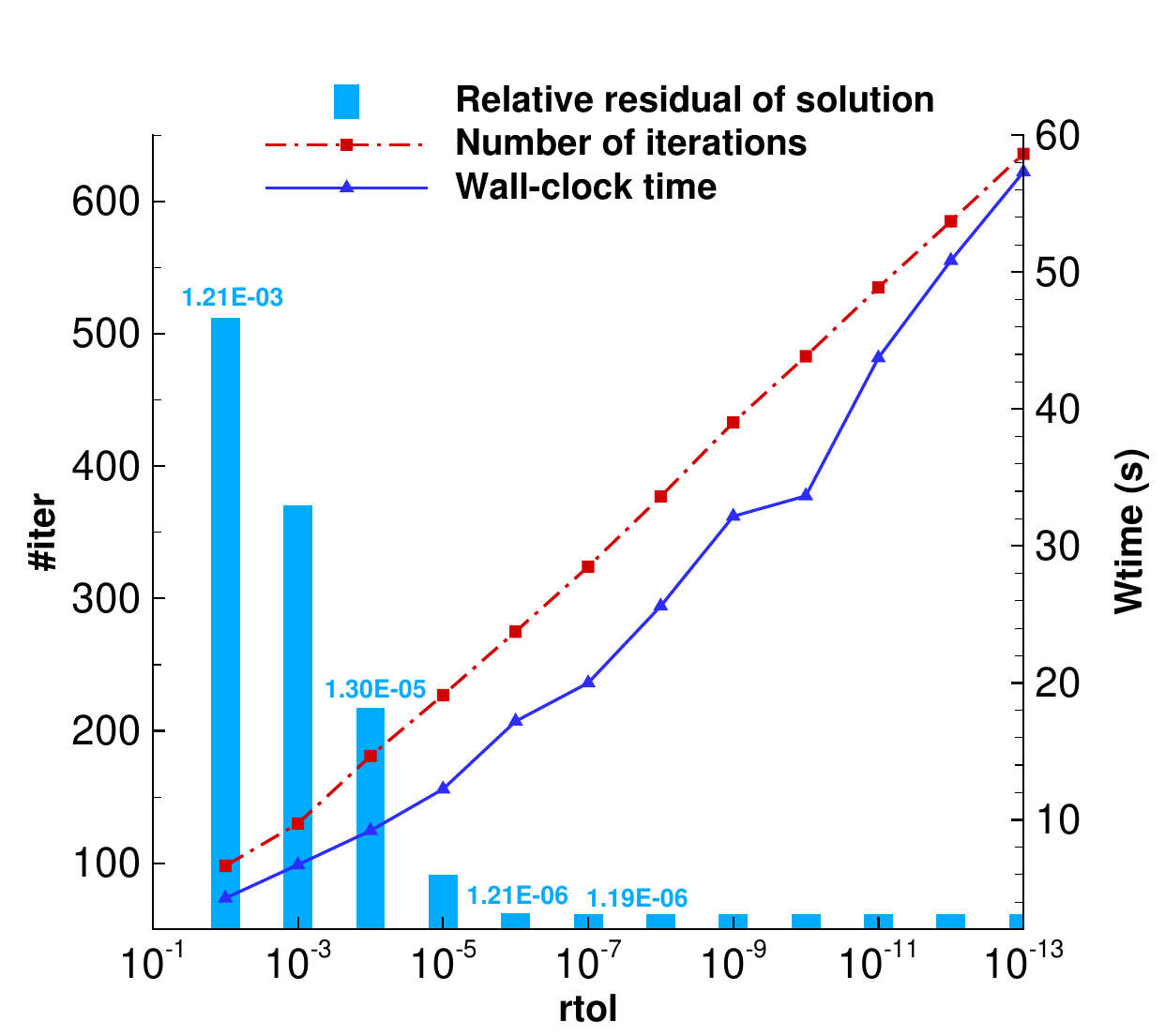}
         \caption{Wedge problem, $f=\SI{40}{\hertz}$, grid size $289 \times 481$.}
         \label{fig:Wedge_on_CoarseSolve_rtol}
     \end{subfigure}
     \hfill
     \begin{subfigure}[b]{0.45\textwidth}
         \centering
         \includegraphics[width=\textwidth]{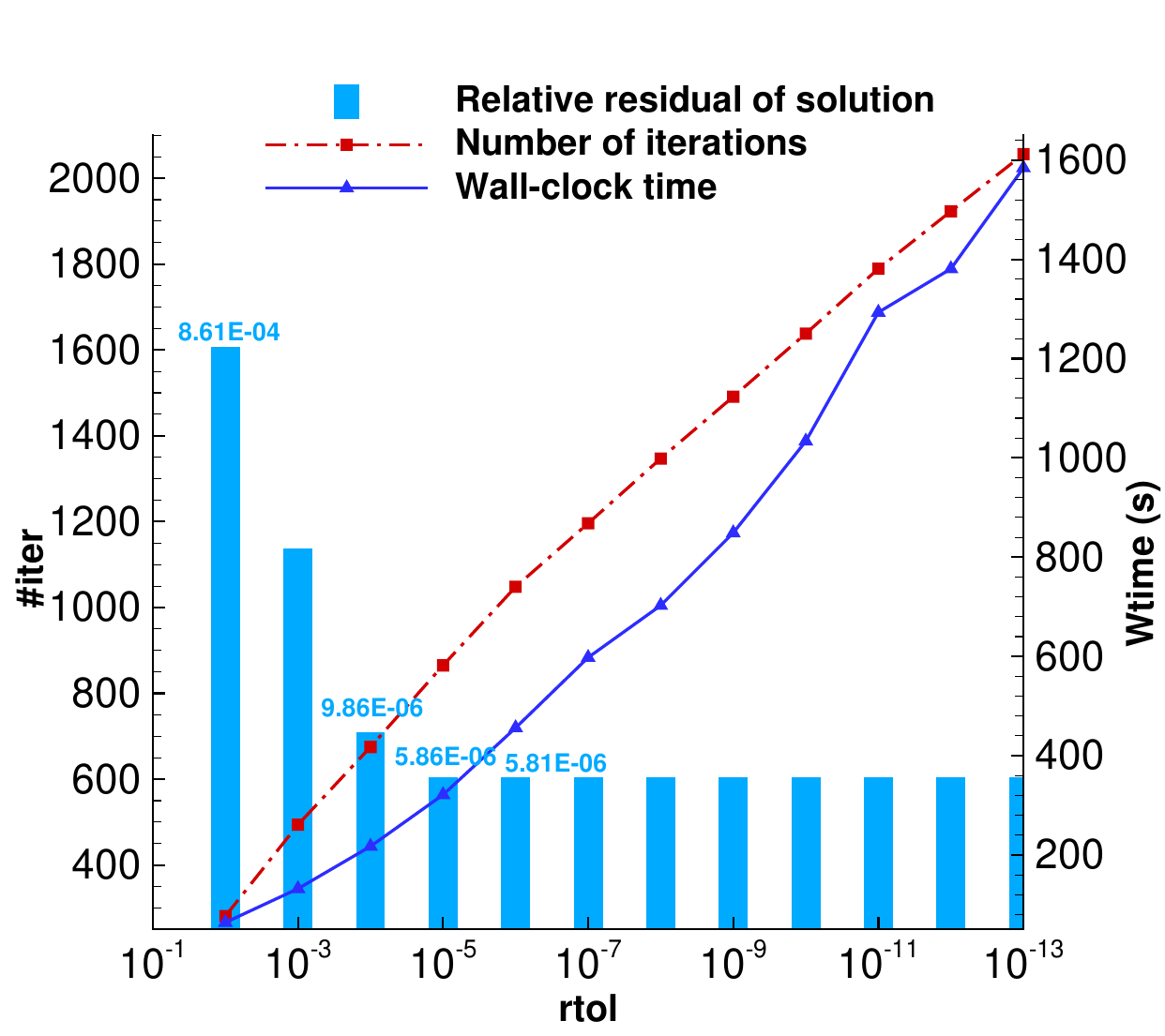}
         \caption{Marmousi problem, $f=\SI{20}{\hertz}$, grid size $1473 \times 481$.}
         \label{fig:Marmousi_on_CoarseSolve_rtol}
     \end{subfigure}
        \caption{Tolerance of the coarse-grid solver for non-constant wavenumber problems. The present APD-preconditioned GMRES is employed for outer iteration.}
    \label{fig:non-constant_on_CoarseSolve_rtol}
\end{figure}

A similar pattern is observed for non-constant wavenumber problems, like the Wedge problem ($f=\SI{40}{\hertz}$, grid size $289 \times 481$) and the Marmousi problem ( $f=\SI{20}{\hertz}$, grid size $1473 \times 481$) as shown in Figure \ref{fig:non-constant_on_CoarseSolve_rtol}.

The findings suggest that, when the present APD-preconditioned GMRES is used, maintaining the tolerance of the coarse-grid solver at the same order of magnitude as that of the outer iterations enables a consistent number of outer iterations and constant relative residual of the final solution while minimizing the total Wall-clock time. It is worth noting that for the GMRES algorithm, we employ left preconditioning and use the preconditioned relative residual as the criterion for convergence. Additional numerical experiments have been conducted to investigate the use of right preconditioning with GMRES, using the unpreconditioned relative residual as the convergence criterion. The results also demonstrate a consistent number of outer iterations across various tolerances for the coarse-grid solver, with only one more iteration compared to the results shown in Table \ref{tab:MP2BC2_outer_iter_varies_CoarseTol}, which strictly reduces the final relative residual to below $10^{-6}$. Nevertheless, we can reach the same conclusion regarding the impact of the coarse-grid solver tolerance on the outer iterations and the relative residual of the solution. We attribute this to that the standard GMRES algorithm for outer iteration, as shown in Algorithm \ref{CSLPDEFGMRES-1}, expects a constant preconditioner. However, if the tolerance of the coarse-grid solver increases beyond a certain threshold, it leads to variable preconditioning. 

The results presented in Figure \ref{fig:Marmousi_on_CoarseSolve_rtol} indicate that solving the coarse-grid problem for the Marmousi model with grid size $1473 \times 481$ and frequency $f=\SI{20}{\hertz}$ requires over 1000 GMRES iterations to achieve a relative residual of $10^{-6}$. To investigate the possibility of using a larger tolerance for the coarse-grid solver and further reduce the solution time, we will proceed with the outer iteration using the GCR algorithm, which allows for a variable preconditioner. While making this transition, all other aspects of the deflation method remain unchanged, except for the use of the GCR algorithm with right preconditioning.

\begin{figure}[h]
    \centering
    \includegraphics[width=0.45\textwidth]{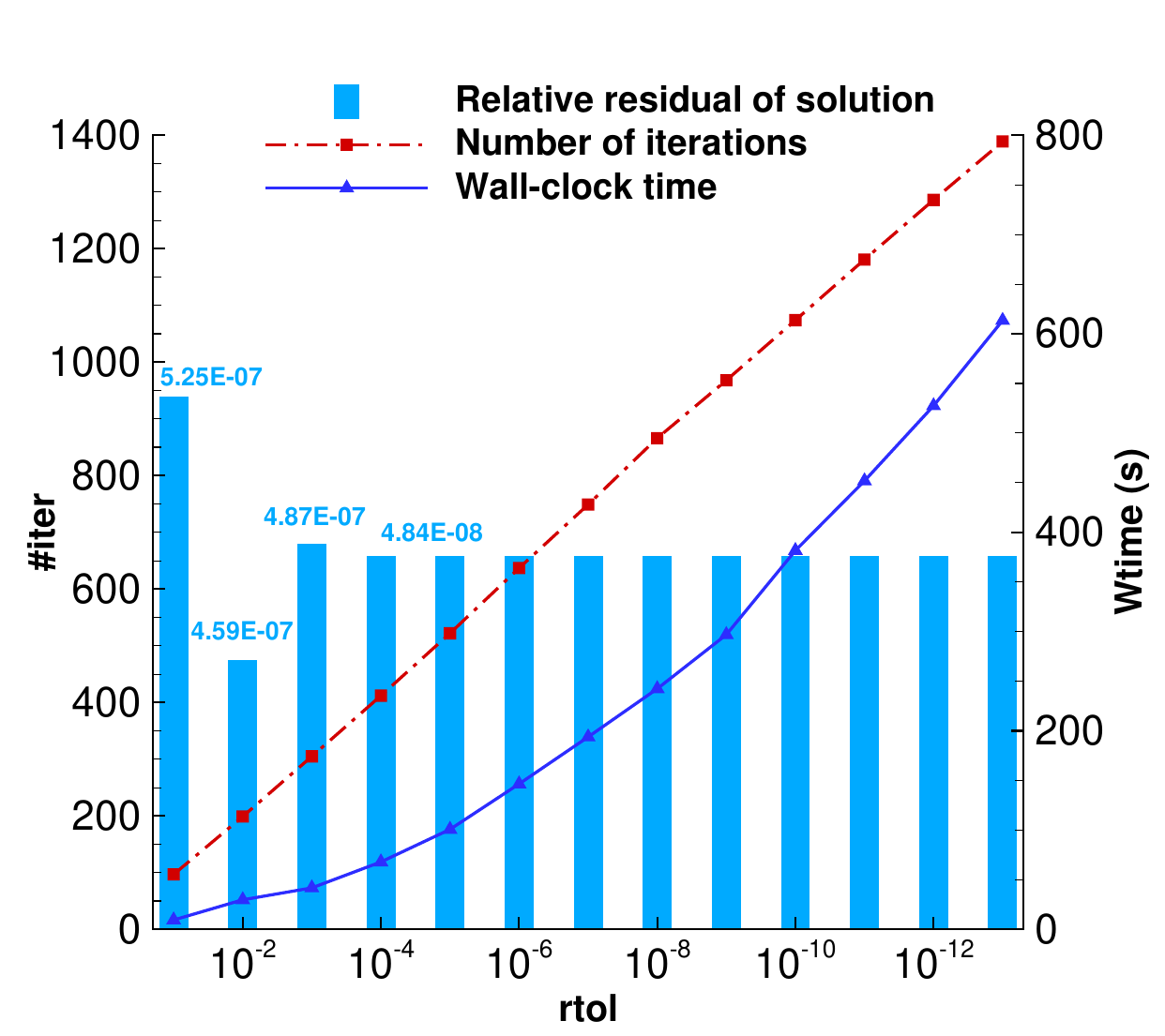}
    \caption{Tolerance of the coarse-grid solver for MP-2b with wavenumber $k=320$ and grid size $513 \times 513$. The present APD-preconditioned GCR is employed for outer iteration.}
    \label{fig:RpGCR_MP2BC2_on_CoarseSolve_rtol}
\end{figure}

For the MP-2b model with wavenumber $k=320$ and grid size $513 \times 513$, the number of outer GCR iterations remains constant at 11, when varying the tolerances of the coarse-grid solver from $10^{-2}$ to $10^{-13}$. If the tolerance of the coarse-grid solver is set to $10^{-1}$, it requires 12 outer iterations. The impact of the coarse-grid solver tolerance on the number of iterations required to solve the coarse-grid problem, the total Wall-clock time, and the relative residual of the final solution is illustrated in Figure \ref{fig:RpGCR_MP2BC2_on_CoarseSolve_rtol}. The variation in the number of iterations for solving the coarse-grid problem and the total Wall-clock time is similar to that shown in Figure \ref{fig:MP2BC2_on_CoarseSolve_rtol}. Similarly, when the tolerance of the coarse-grid solver is smaller than the tolerance of the outer iterations, the residual of the solution remains unchanged. However, when the tolerance of the coarse-grid solver exceeds the tolerance of the outer iterations, the residual of the solution remains within the tolerance instead of increasing, as observed in the standard GMRES method. The same behavior is observed in the case of the non-constant wavenumber model problems. Thus, using the GCR algorithm as the outer iteration method allows us to set the tolerance of the coarse-grid solver to a relatively large value $10^{-1}$, while still maintaining the desired accuracy of the final solution. For this MP-2b model problem, if the APD-preconditioned GMRES is employed and the tolerance for the coarse-grid solver is set to $10^{-6}$, the total Wall-clock time is $\SI{145.41}{\second}$.  If the APD-preconditioned GCR is employed and the tolerance for the coarse-grid solver is set to $10^{-1}$, although one extra outer iteration is required, the total Wall-clock time is $\SI{9.48}{\second}$.  The total computation time is reduced by around 95\% of the time required when using a tolerance of $10^{-6}$ for the coarse-grid solver.  

Additional numerical experiments have provided further evidence that regardless of whether the standard GMRES algorithm is employed with left or right preconditioning as the outer iteration algorithm, if the tolerance of the coarse-grid solver is set larger than that of the outer iterations, the residual of the final solution will increase. However, when utilizing the ADP-preconditioned Flexible GMRES algorithm, which also allows for variable preconditioning, similar to the ADP-preconditioned GCR approach, even with a relatively larger tolerance for the coarse-grid solver, the accuracy of the final solution can still be maintained within the specified tolerance range.

In the subsequent parallel performance study of the parallel APD-preconditioned GMRES and GCR, we will set the tolerance for the coarse-grid problem solver to $10^{-6}$ and $10^{-1}$, respectively. 
% Though additional numerical experiments suggest that a stricter tolerance for the coarse-grid problem solver does not significantly impair the parallel performance.

\subsection{Parallel performance}
In this section, we conduct a comprehensive evaluation of the parallel performance of our solver, focusing on weak scaling and strong scaling. Weak scaling refers to the solver's performance when the problem size and the number of processing elements increase proportionally. This allows us to assess our parallel solver's ability to solve large-scale heterogeneous Helmholtz problems with minimized pollution error by using smaller step size $h$. On the other hand, strong scaling measures the capability of our parallel solution method to solve a fixed problem size more quickly by adding more resources. Through these scaling assessments, we aim to understand the effectiveness of our solver in handling various problem sizes and resource allocations, providing crucial insights into its optimal operational parameters and potential scalability limits.

\subsubsection{Weak scaling}
\begin{figure}
     \centering
     \begin{subfigure}[b]{0.45\textwidth}
         \centering
         \includegraphics[width=\textwidth]{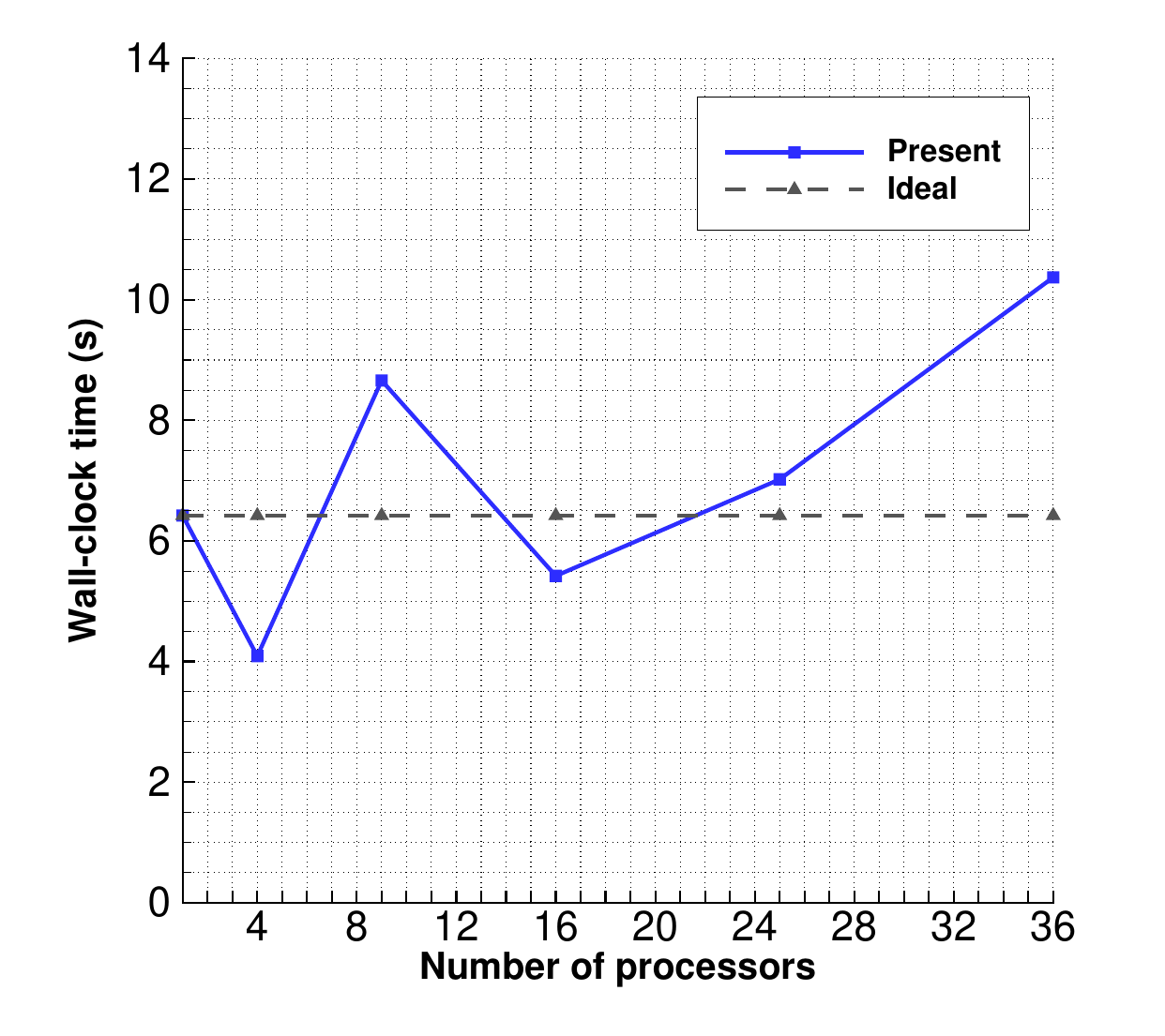}
         \caption{APD-preconditioned GMRES}
         \label{fig:weak-scaling-MP2bk100-gmres}
     \end{subfigure}
     \hfill
     \begin{subfigure}[b]{0.45\textwidth}
         \centering
         \includegraphics[width=\textwidth]{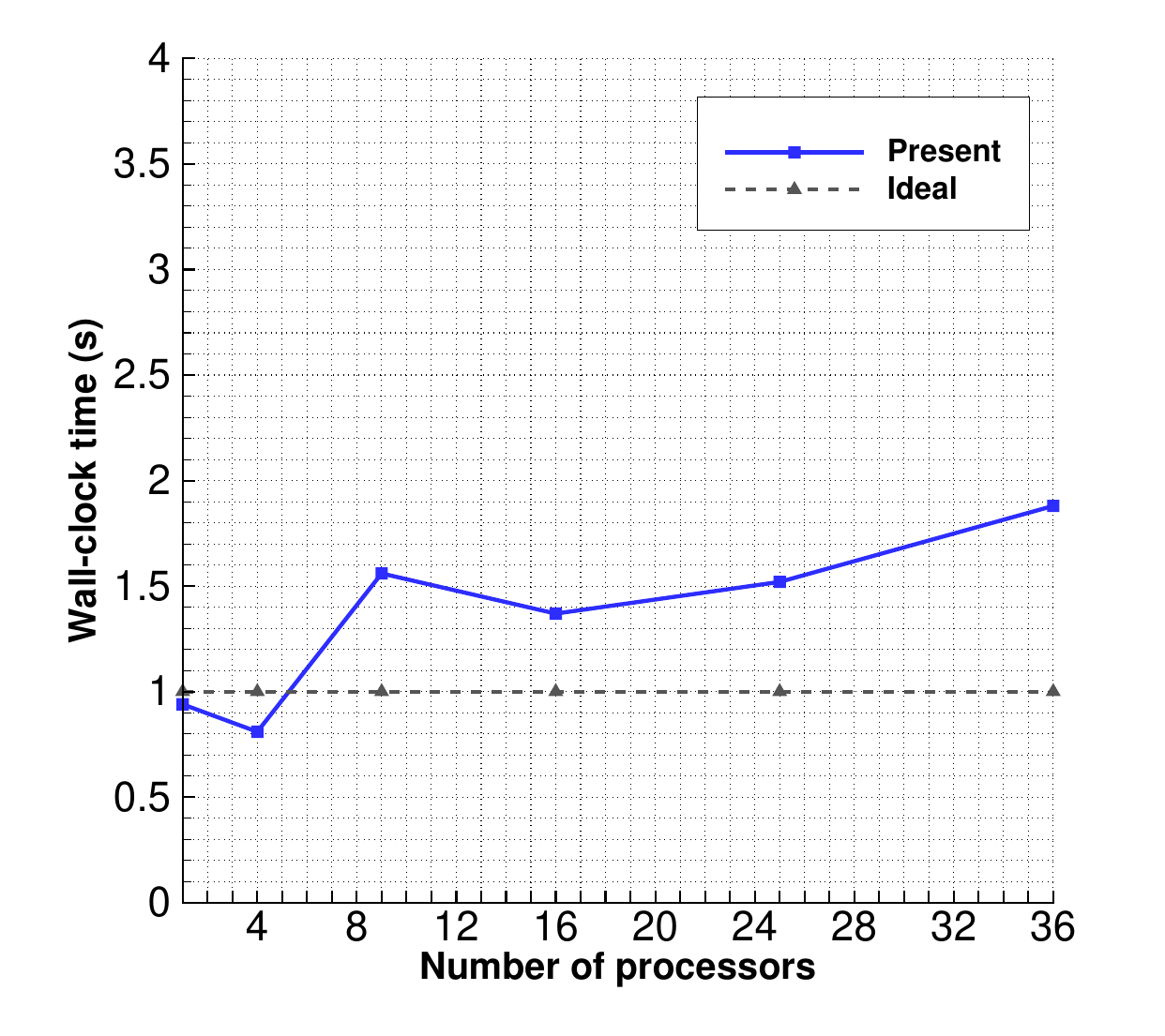}
         \caption{APD-preconditioned GCR}
         \label{fig:weak-scaling-MP2bk100-gcr}
     \end{subfigure}
        \caption{Weak scaling for MP-2b with $k=100$ and a grid size of $160 \times 160$ per processes.}
    \label{fig:weak-scaling-MP2bk100}
\end{figure}

Figure \ref{fig:weak-scaling-MP2bk100} shows the results for the weak scalability test solving the MP-2b model problem with $k=100$ from 1, 4, 9, ..., up to 36 processes. The problem size was refined from $161 \times 161$, $321 \times 321$, $481 \times 481$, ... up to $961 \times 961$, ensuring that each process handled a grid size of approximately $160 \times 160$. Clearly, for the same problem and with the same number of processes, the ADP-preconditioned GCR method with a coarse-grid solver tolerance of $10^{-1}$ exhibits a significantly reduced computational time compared to the ADP-preconditioned GMRES method with a coarse-grid solver tolerance of $10^{-6}$.

From Figure \ref{fig:weak-scaling-MP2bk100}, it appears that as the grid size and the processes increase proportionally ($np \leq 25$), the required Wall-clock time does not remain perfectly constant. This behavior may be attributed to the relatively small size of the model problem and the short computation time ($\approx \SI{1}{\second}$ for GCR), resulting in a significant proportion of communication time in the total elapsed time.
Additionally, as this numerical experiment was conducted within a single compute node with 48 cores, the Wall-clock time increases as the utilization of cores approaches 48, primarily due to the increased data communication and limited bandwidth within a single compute node. The same trend holds for model problems with non-constant wavenumbers, as shown in Table \ref{tab:non-constant-weak-scaling-gmres} and \ref{tab:non-constant-weak-scaling-gcr}. The present weak scalability meets our requirements for minimizing pollution error by grid refinement within a certain range.

\begin{table}[htb]
\centering
\caption{Weak scaling for model problems with non-constant wavenumber. The present APD-preconditioned GMRES is employed for outer iteration.}
\scalebox{0.85}{
\begin{tabular}{cccc}
\hline
grid size & np & \#iter   & CPU time (s) \\ \hline
\multicolumn{4}{c}{Wedge, $f=\SI{40}{\hertz}$}           \\
577  $\times$ 961   & 6  & 9 (251)  & 53.41        \\
1153 $\times$ 1921  & 24 & 10 (259) & 68.59        \\
          &    &          &              \\
\multicolumn{4}{c}{Marmousi, $f=\SI{10}{\hertz}$}        \\
737  $\times$ 241  & 3  & 10 (414) & 103.92       \\
1473 $\times$ 481  & 12 & 9 (378)  & 111.20       \\
2945 $\times$ 961  & 48 & 9 (383)  & 156.45       \\ \hline
\end{tabular}}
\label{tab:non-constant-weak-scaling-gmres}
\end{table}

\begin{table}[htbp]
\centering
\caption{Weak scaling for model problems with non-constant wavenumber. The present APD-preconditioned GCR is employed for outer iteration.}
\scalebox{0.85}{
\begin{tabular}{cccc}
\hline
grid size & np & \#iter   & CPU time (s) \\ \hline
\multicolumn{4}{c}{Wedge, $f=\SI{40}{\hertz}$}           \\
577  $\times$ 961   & 6  & 10 (46)  & 4.86        \\
1153 $\times$ 1921  & 24 & 10 (43)  & 5.75        \\
          &    &          &              \\
\multicolumn{4}{c}{Marmousi, $f=\SI{10}{\hertz}$}        \\
737  $\times$ 241  & 3  & 11 (63) & 10.55       \\
1473 $\times$ 481  & 12 & 10 (58)  & 12.08       \\
2945 $\times$ 961  & 48 & 10 (58)  & 17.72       \\ \hline
\end{tabular}}
\label{tab:non-constant-weak-scaling-gcr}
\end{table}

\subsubsection{Strong scaling}
\begin{figure}
     \centering
     \begin{subfigure}[h]{0.45\textwidth}
         \centering
         \includegraphics[width=\textwidth]{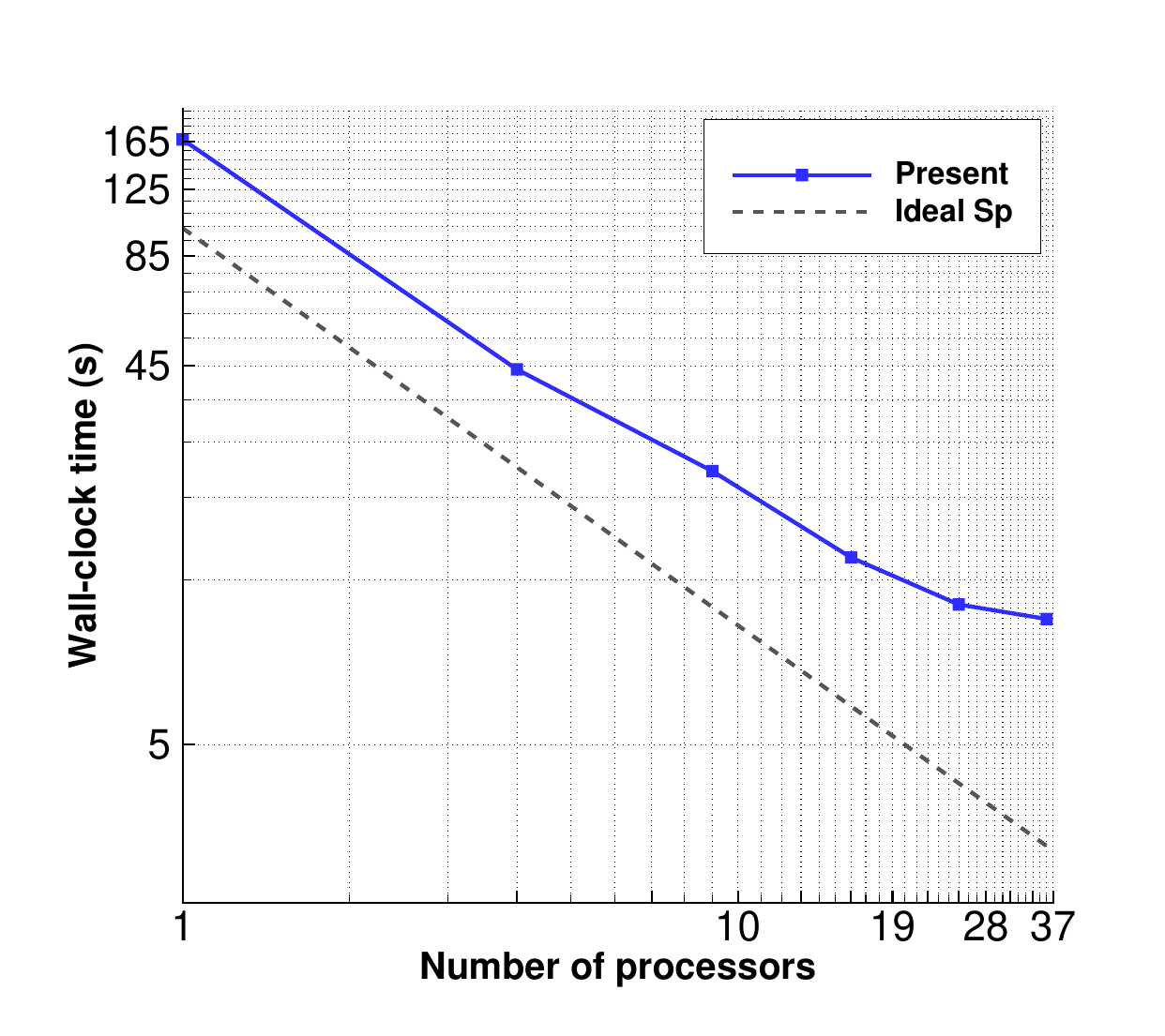}
         \caption{APD-preconditioned GMRES}
         \label{fig:strong-scaling-MP2bk100-gmres}
     \end{subfigure}
     \hfill
     \begin{subfigure}[h]{0.45\textwidth}
         \centering
         \includegraphics[width=\textwidth]{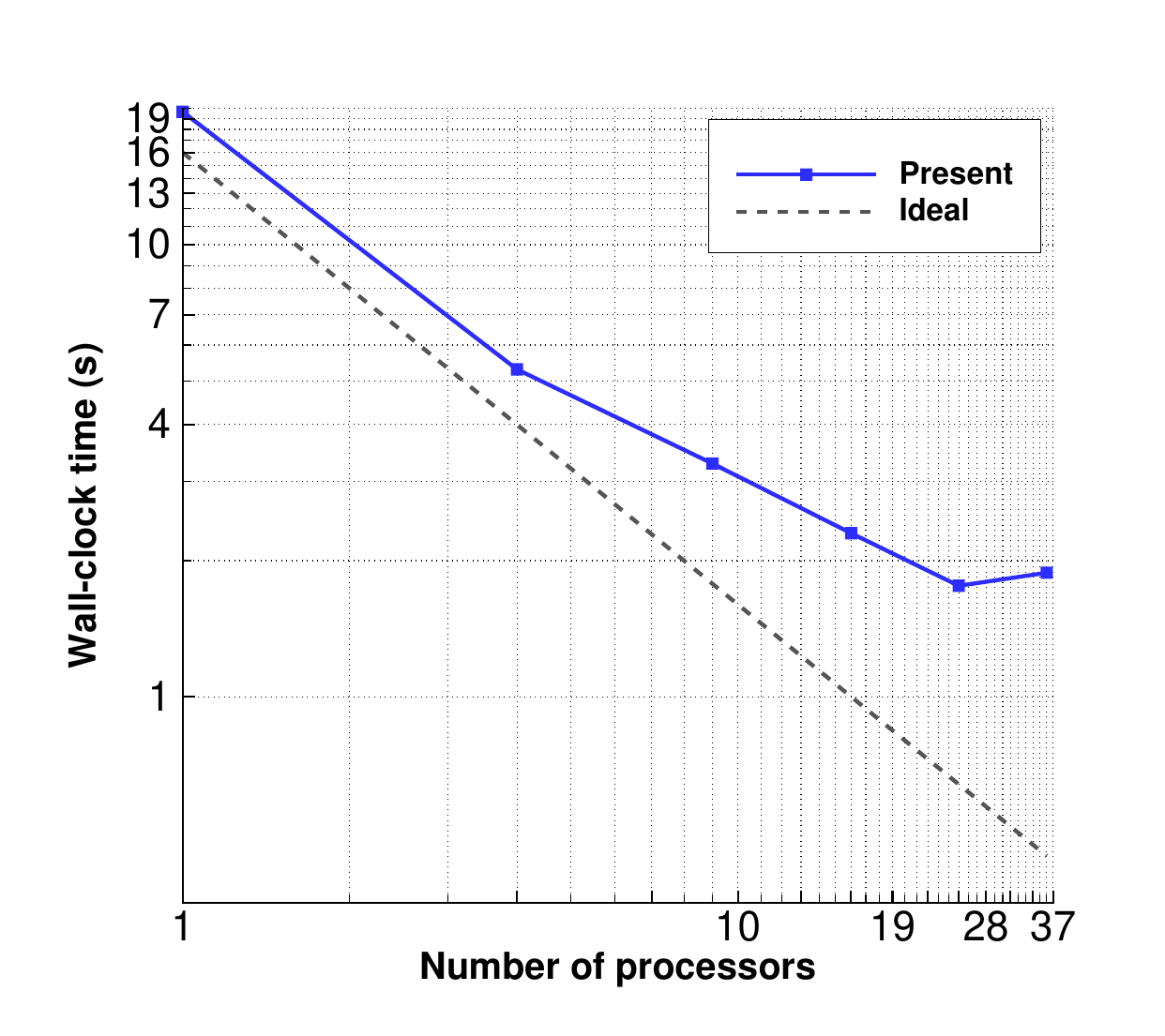}
         \caption{APD-preconditioned GCR}
         \label{fig:strong-scaling-MP2bk100-gcr}
     \end{subfigure}
        \caption{Strong scaling for MP-2b with $k=100$ and a grid size of $961 \times 961$.}
    \label{fig:strong-scaling-MP2bk100}
\end{figure}

\begin{figure}
     \centering
     \begin{subfigure}[h]{0.45\textwidth}
         \centering
         \includegraphics[width=\textwidth]{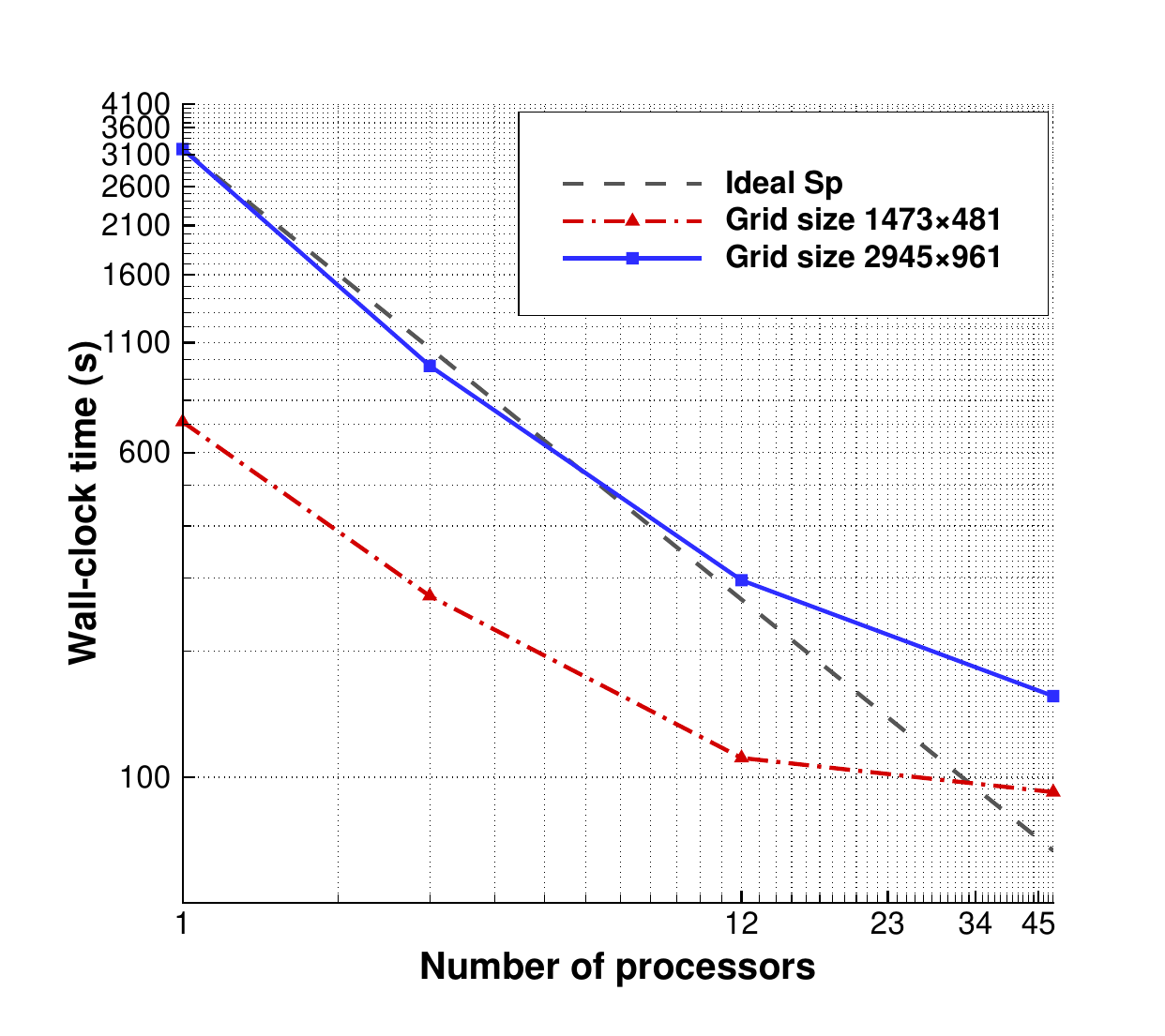}
         \caption{APD-preconditioned GMRES}
         \label{fig:strong-scaling-Marmousi-gmres}
     \end{subfigure}
     \hfill
     \begin{subfigure}[h]{0.45\textwidth}
         \centering
         \includegraphics[width=\textwidth]{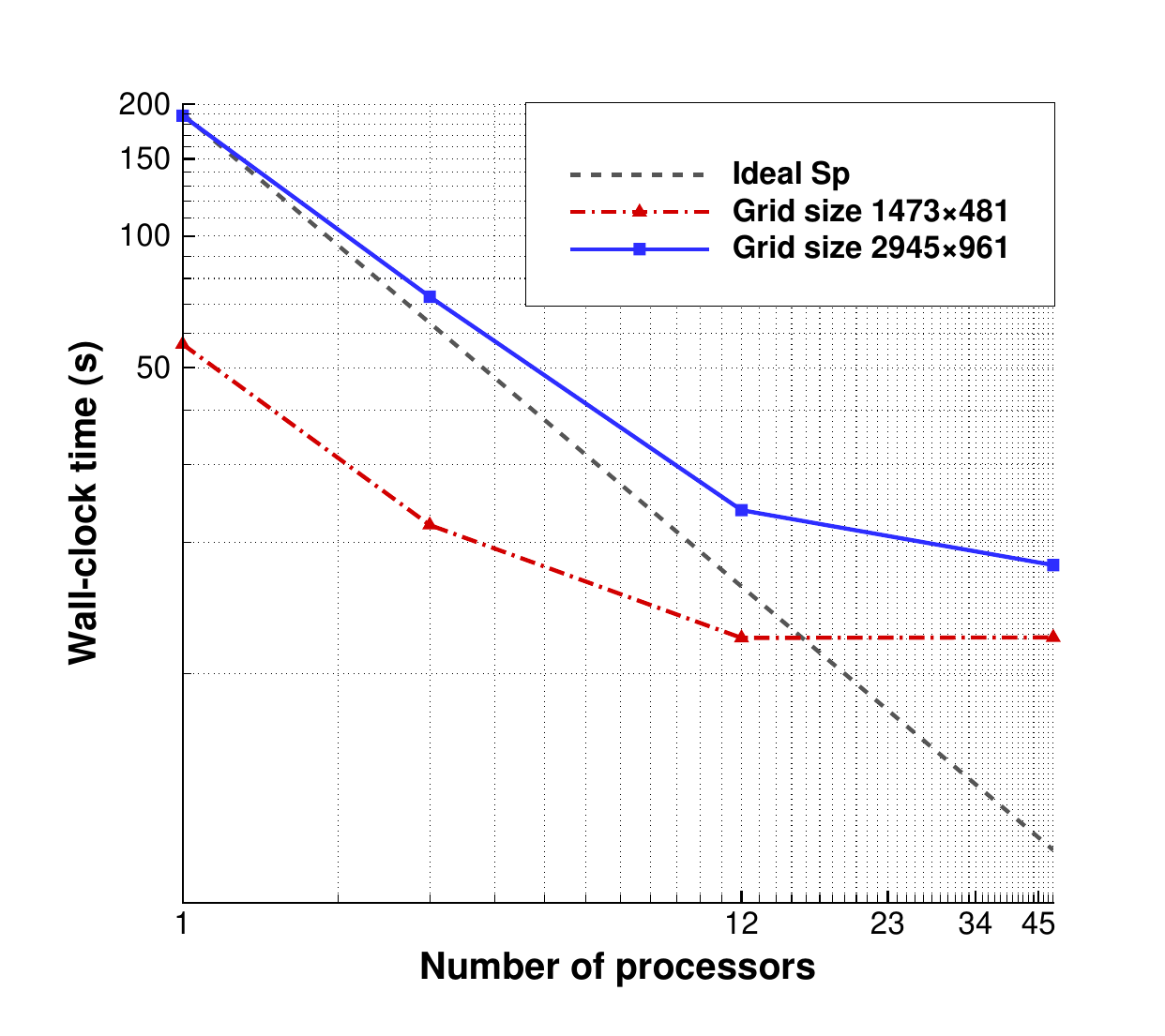}
         \caption{APD-preconditioned GCR}
         \label{fig:strong-scaling-Marmousi-gcr}
     \end{subfigure}
        \caption{Strong scaling for Marmousi problem with $f=\SI{10}{\hertz}$.}
    \label{fig:strong-scaling-Marmousi}
\end{figure}

We are also interested in the strong scalability properties of the present parallel deflation method for the Helmholtz problems. First of all, numerical experiments show that the number of external iterations required is found to be independent of the number of processes, which is a favorable property of our solution method. We consider the MP-2b problem with a wavenumber of $k=100$ and a grid size of $961 \times 961$ on a growing number of processes. Figure \ref{fig:strong-scaling-MP2bk100} presents the Wall-clock time versus the number of processes. The numerical experiment conducted on close-to-48 processes exhibits a moderate decrease in terms of parallel efficiency. This can be partly attributed to the increased amount of communication, resulting in a significant decrease in the computation/communication ratio. Similar patterns can be observed for the Marmousi problem with a frequency of $f=\SI{10}{\hertz}$ and a grid size of $1473 \times 481$, as shown by the red dotted line in Figure \ref{fig:strong-scaling-Marmousi}. When solving larger grid sizes, such as $2945 \times 961$ as indicated by the solid blue line in Figure \ref{fig:strong-scaling-Marmousi}, the computation/communication ratio increases, leading to improved parallel efficiency compared to the case with a grid size of $1473 \times 481$.

\begin{figure}
     \centering
     \begin{subfigure}[h]{0.45\textwidth}
         \centering
         \includegraphics[width=\textwidth]{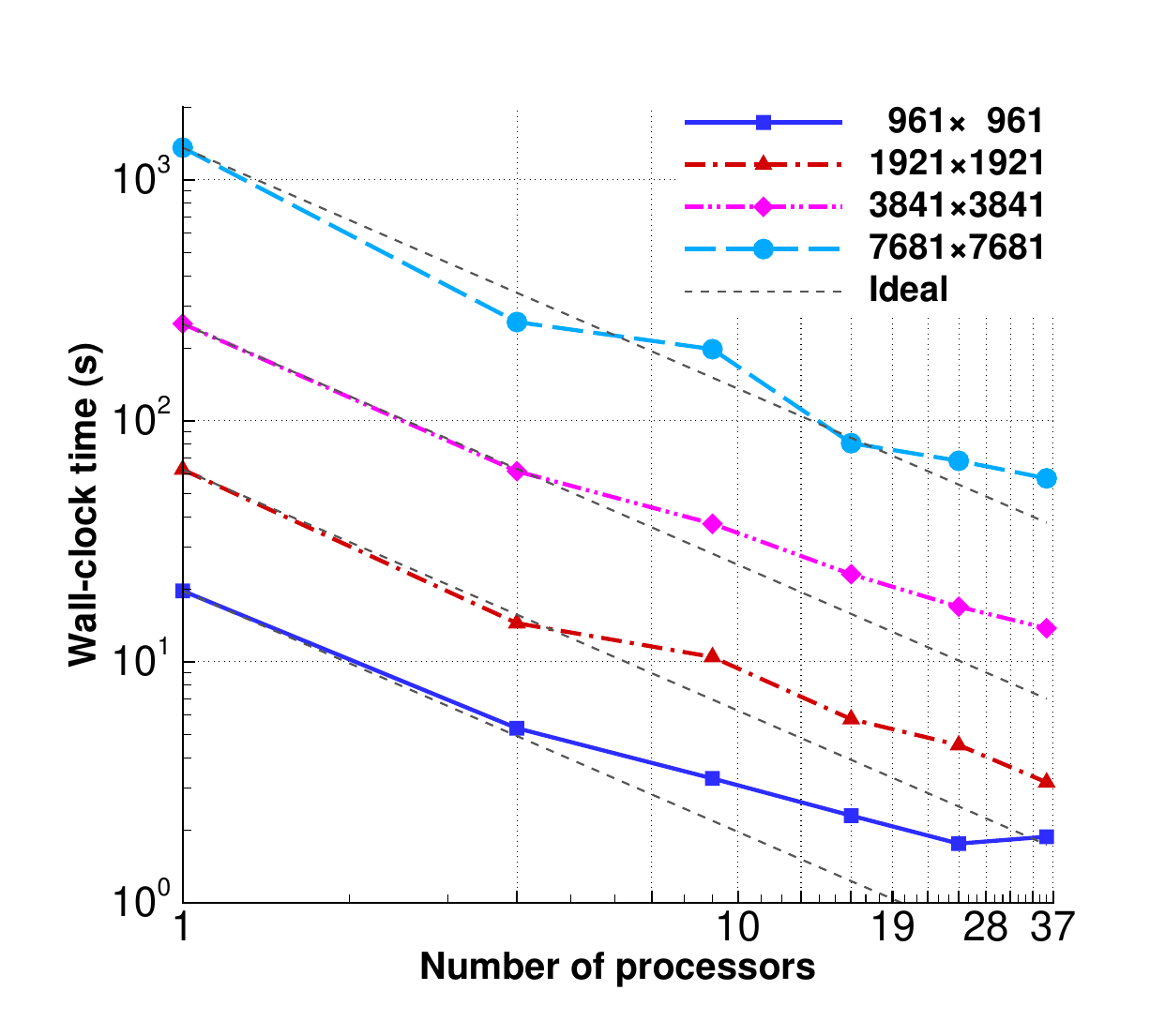}
         \caption{$k=100$}
         \label{fig:strong-scaling-MP2bk100-gcr_girds}
     \end{subfigure}
     \hfill
     \begin{subfigure}[h]{0.45\textwidth}
         \centering
         \includegraphics[width=\textwidth]{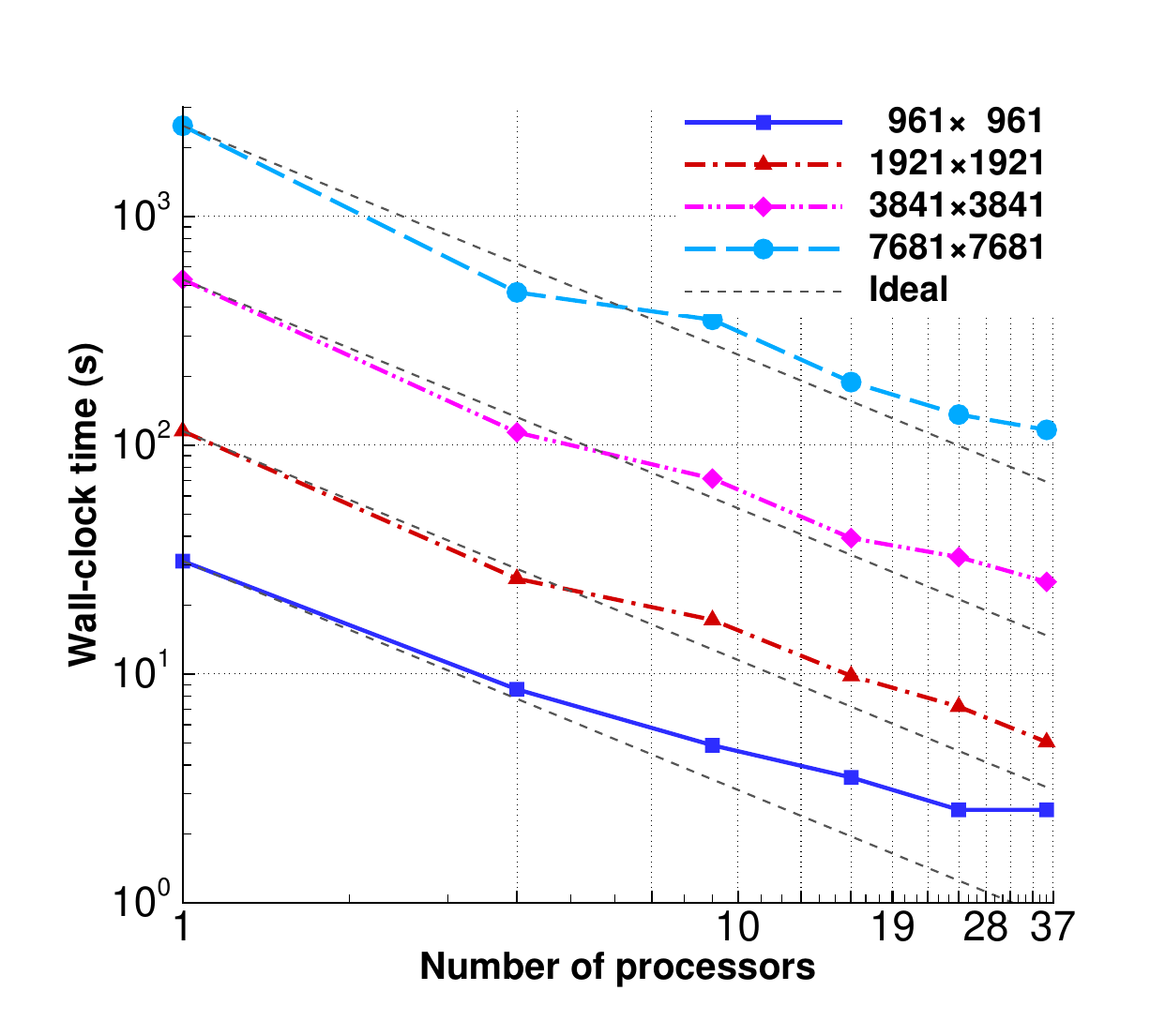}
         \caption{$k=200$}
         \label{fig:strong-scaling-MP2bk200-gcr_girds}
     \end{subfigure}
        \caption{Strong scaling of APD-preconditioned GCR for MP-2b with various grid sizes.}
    \label{fig:strong-scaling-MP2bk100k200-gcr_girds}
\end{figure}

One can find that the strong scaling of the present parallel APD-preconditioned GMRES and GCR behave in a similar way. In terms of parallel efficiency, using GMRES as outer iteration performs a bit better than using GCR. It is not because of the inherent parallelism of the algorithm but the relatively smaller computation/communication ratio for the GCR scenario. If we increase the grid size, see Figure \ref{fig:strong-scaling-MP2bk100k200-gcr_girds}, an improvement in parallel efficiency can be observed. This is because the number of ghost-grid layers used for communication remains constant, the amount of data to be communicated doubles when the number of grid points doubles in each direction. However, the total number of grid points increases fourfold, resulting in a larger ratio of computational and communication time. Furthermore, we also observe that the computational time increases nearly proportionally to the grid size, approximately by a factor of $4$-$5$. But the impact of communication overhead becomes increasingly evident. From Figure \ref{fig:strong-scaling-MP2bk100k200-gcr_girds}, one can find that similar strong scalability holds for a larger wavenumber.

\begin{figure}
     \centering
     \begin{subfigure}[h]{0.45\textwidth}
         \centering
         \includegraphics[width=\textwidth]{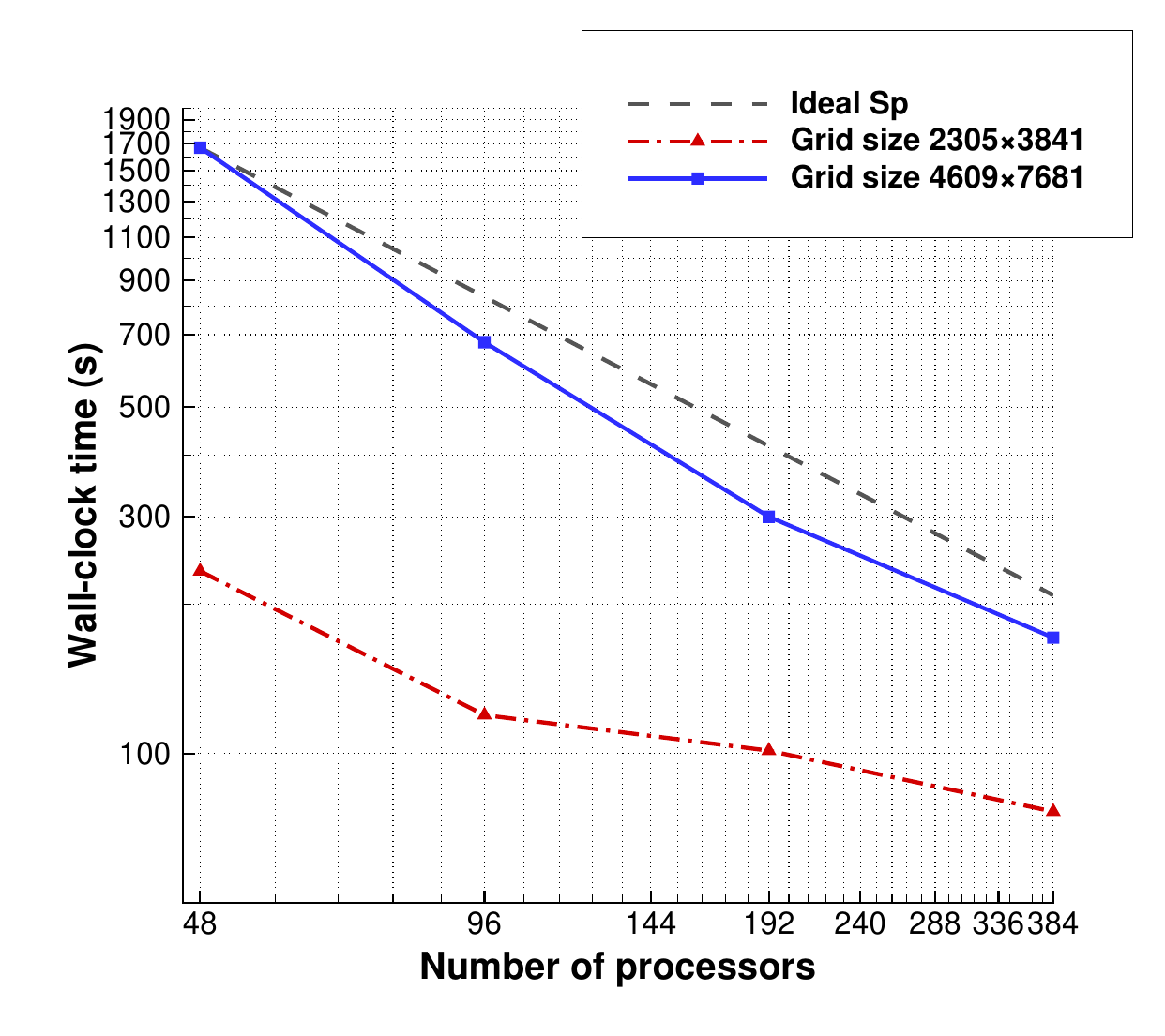}
         \caption{APD-preconditioned GMRES}
         \label{fig:strong-scaling-Wedge-gmres}
     \end{subfigure}
     \hfill
     \begin{subfigure}[h]{0.45\textwidth}
         \centering
         \includegraphics[width=\textwidth]{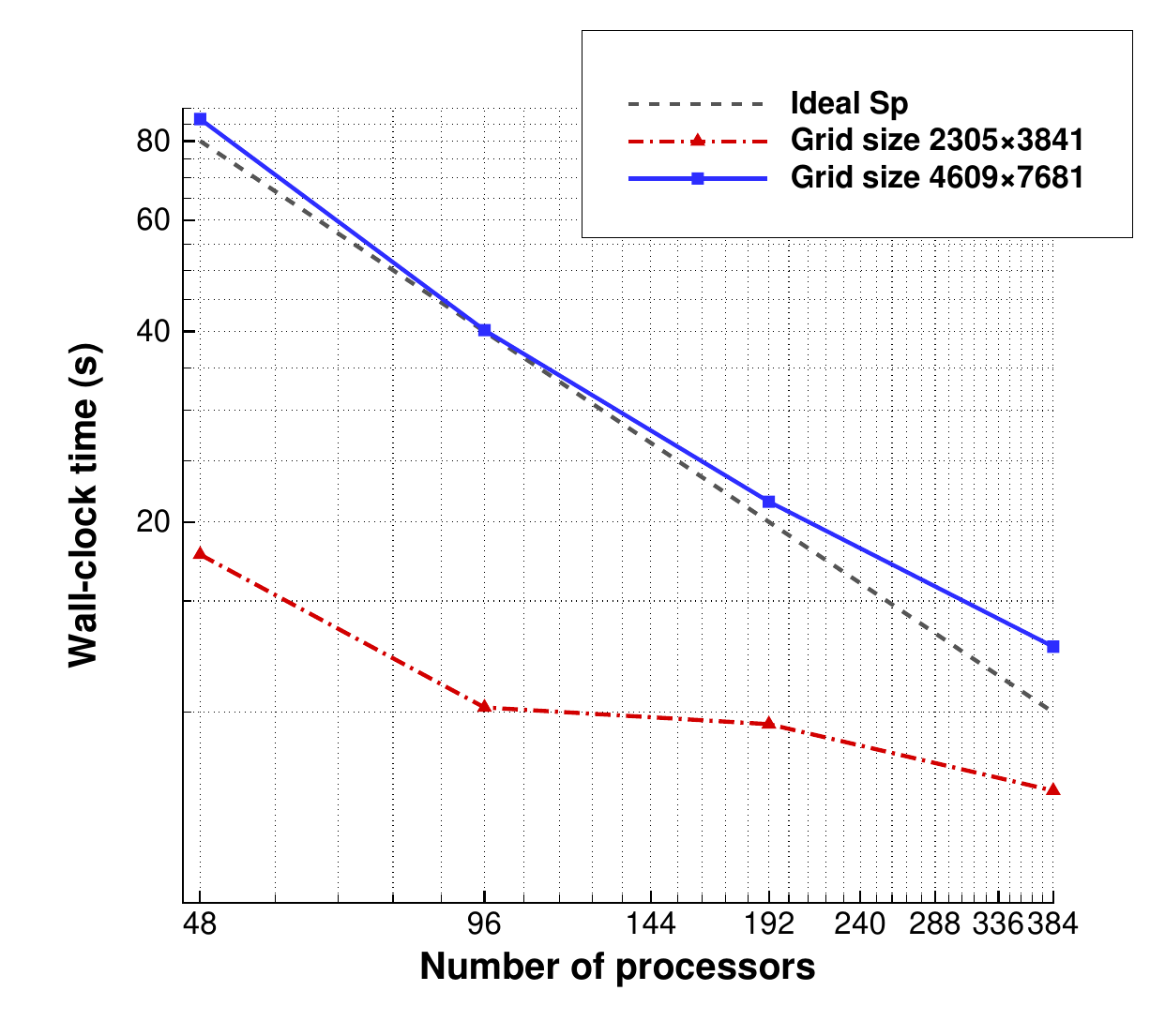}
         \caption{APD-preconditioned GCR}
         \label{fig:strong-scaling-Wedge-gcr}
     \end{subfigure}
        \caption{Strong scaling for Wedge problem with $f=\SI{40}{\hertz}$.}
    \label{fig:strong-scaling-Wedge}
\end{figure}

\begin{figure}[h]
    \centering
    \includegraphics[width=0.45\textwidth]{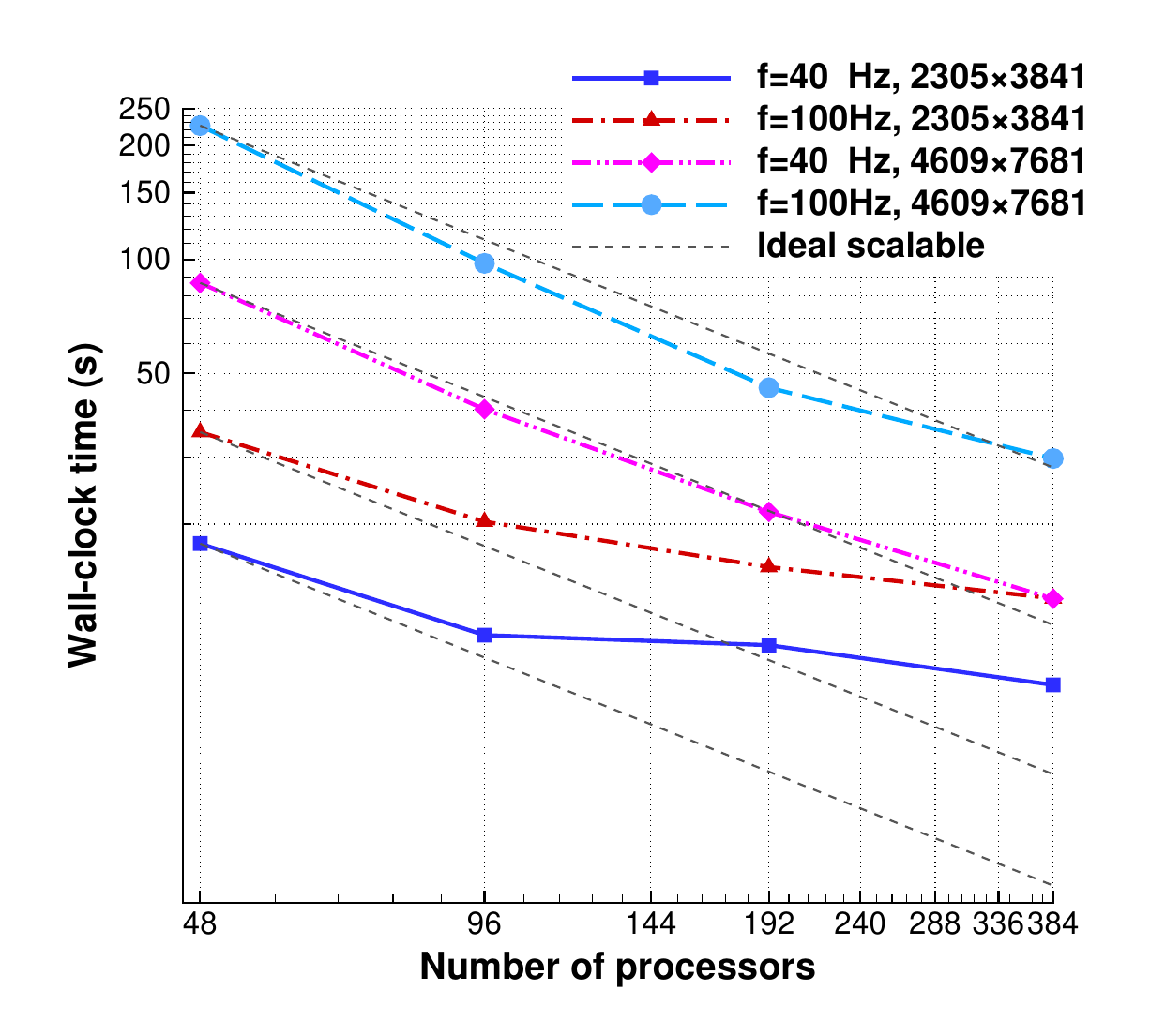}
    \caption{Strong scaling of APD-preconditioned GCR for Wedge problem with $f=\SI{40}{\hertz}$ and $f=\SI{100}{\hertz}$.}
    \label{fig:strong-scaling-wedgef40vs100-gcr_girds}
\end{figure}

The strong scalability property of the present solver on multiple compute nodes is also investigated. Figure \ref{fig:strong-scaling-Wedge} illustrates the time required to solve the wedge model problem with a frequency of $f=\SI{40}{\hertz}$ on at most 6 compute nodes (a total of 384 processes). For a grid size of $2305 \times 3841$, the numerical experiment performed on more than $96$ processes exhibits a decrease in parallel efficiency. However, as the grid size increases to $4609 \times 7681$, which consequently increases the computation/communication ratio, the parallel efficiency significantly improves, with speedup even surpassing the ideal case, i.e., superlinear speedup. When parallel computing tasks can fully utilize the caches on multiple compute nodes, data access speeds are faster, thus enhancing computational efficiency. 
For solving the Wedge problem with a higher frequency of $f=\SI{100}{\hertz}$, as depicted in Figure \ref{fig:strong-scaling-wedgef40vs100-gcr_girds}, we observe similar patterns to $f=\SI{40}{\hertz}$. Despite the constant number of outer iterations, the number of iterations required to solve the coarse-grid problem increases from around 45 to around 110. Consequently, the computation time also increases by a factor of 2-3. In terms of parallel efficiency, we observe similar behavior.

\section{Conclusions}
\label{sec:conclusions}
We propose a matrix-free parallel two-level deflation preconditioner for solving two-dimensional Helmholtz problems in both homogeneous and heterogeneous media.

By leveraging the matrix-free parallel framework and geometric multigrid-based CSLP \cite{jchen2D2022}, we provide a matrix-free parallelization method for the general two-level deflation preconditioner for the Helmholtz equation discretized with finite differences. Numerical experiments demonstrate that, compared with the Galerkin coarsening approach, the method of using re-discretization to obtain the coarse grid operator of deflation will slow down convergence in a certain.

To enhance the convergence properties, the higher-order approximation scheme proposed by Dwarka and Vuik \cite{dwarka2020scalable} is employed to construct the deflation vectors. The study presents higher-order deflation vectors in two dimensions and their matrix-free parallel implementation. Various methods for implementing matrix-free coarse-grid operators are proposed and compared. The numerical experiments demonstrate the effectiveness of the proposed coarse-grid re-discretization scheme based on the Galerkin coarsening approach, which achieves wave-number-independent convergence for both constant and non-constant wave number model problems.

Furthermore, we have studied in detail the tolerance setting of the coarse grid solver when the deflation preconditioned GMRES type methods are employed to solve the Helmholtz problems. An optimal tolerance for the coarse grid solver can effectively reduce the solution time. 

Finally, the performance of the present parallel GCR solver preconditioned by higher-order deflation is studied, including both weak and strong scalability. Using a matrix-free approach to reduce memory requirements and weak scalability allows us to minimize pollution errors by refining the grid. Strong scalability allows us to solve high-frequency heterogeneous Helmholtz problems faster and more effectively in parallel. The presented two-level preconditioner serves as a benchmark for potential future developments of matrix-free parallel multilevel deflation preconditioners.

\clearpage
\appendix

\section{A potential way to construct a nine-point compact re-discretization scheme}  \label{Appendix B}
From ReD-Glk, we know that the re-discretization on the coarse grid does not have to be high order accuracy. To guarantee that formulation \cref{eq:GeneralFormReD} is at least a low-order approximation of the Helmholtz operator, the coefficients should satisfy that
\begin{eqnarray}
    \label{eq:CommonConstraits_A}
    &&a_0+4a_s+4a_c=0, \quad a_s+2a_c=-1, \\
    \label{eq:CommonConstraits_B} &&b_0+4b_s+4b_c=1.
\end{eqnarray}

 \cite{dwarka2020scalable} reveals that to obtain a close-to-wavenumber-independent convergence, it is crucial that the near-zero eigenvalues of the coarse-grid operators remain aligned with those of the fine-grid operator. For the fine-grid operator, if we use the second-order finite difference scheme \cref{eq:Helmstencil} and homogeneous Dirichlet conditions, the discrete eigenvalues are given by 
\begin{equation}
    \hat{\lambda}_h^{i,j}=\frac{1}{h^2}\left[ 4-2\cos(i\pi h)-2\cos(j \pi h)-k^2 h^2\right], \quad i,j=1,2,...,n-1.
\end{equation}
If we use formulation \cref{eq:GeneralFormReD} and homogeneous Dirichlet conditions, one can also derive the analytical eigenvalues as following
\begin{equation}
\label{eq:eigs_A2h_optstcl}
    \begin{aligned}
    \hat{\lambda}_{H}^{p,q}&=\frac{1}{H^2}\left[ (a_0-b_0k^2H^2)+2(a_s-b_sk^2H^2)(\cos(p\pi H)+\cos(q \pi H)) \right.\\
    &\left. +4(a_c-b_ck^2H^2)\cos(p\pi H)\cos(q \pi H)\right], \quad p,q=1,2,...,n/2-1.
    \end{aligned}
\end{equation}
where $H=2h$. 

Together with the constraints above, we can have an optimization problem to determine the coefficients in \cref{eq:GeneralFormReD}. That is, for certain large wavenumber $k$, the near-zero value of $\hat{\lambda}_{H}^{p,q}$ remains aligned with the near-zero value of $\hat{\lambda}_h^{i,j}$.

Let us take $k=80$ for example. Suppose $h=0.00001$, we can have
	\begin{equation*}
		\lim_{h \rightarrow 0} \lambda_h^{i,j} \approx \hat{\lambda}_h^{i,j}
	\end{equation*}
After sorting in ascending order, \[\min |\lambda_h^{i,j}|= |-4.496|\] corresponding to index $482$, and \[\min|\hat{\lambda}_h^{i,j}|=|-4.513|\] also corresponding to index $482$. From \[(i^2+j^2)\pi^2-k^2=-4.496,\] we can get that $i=18$ and $j=18$. One can verify that $\hat{\lambda}_h^{18,18} = -4.513$. Thus, let $p=18$ and $q=18$ for \cref{eq:eigs_A2h_optstcl}, we would like that \[\hat{\lambda}_{H}^{18,18} = -4.5\]
	Combined with Eq. (\ref{eq:CommonConstraits_A}) and (\ref{eq:CommonConstraits_B}), we have only four constrains for six unknown coefficients. For simplicity, let $b_s=b_c=0$ and $b_0=1$. Then we can solve that 
	\begin{equation}
		\label{eq:optstcl_k80}
		a_0=4.632 \quad a_s=-1.316 \quad a_c=0.158
	\end{equation} 
Let us denote this stencil as \textbf{ReD-9ptO2}. With the given coefficients, for $k=100$, one can get
	\[\min |\lambda_h^{i,j}|\approx |-2.1|,\quad \text{index}=\{763, 764\}, \quad i=22, j=23\]
	\[\min|\hat{\lambda}_h^{i,j}|\approx |-2.1|,\quad \text{index}=\{763, 764\}, \quad i=22, j=23\]
	\[\min|\hat{\lambda}_{H}^{p,q}| \approx  |-2.1|,\quad \text{index}=\{763, 764\}, \quad p=22, q=23\]

	However, in practice, $h$ will not be so small. For $k=80$ and $kh=0.3125$, that is $h=\frac{1}{256}$, we will have
	\[\min |\lambda_h^{i,j}| = |-4.496|,\quad \text{index}=\{482\}\]
	\[\min|\hat{\lambda}_h^{i,j}| = |-11.200|,\quad \text{index}=\{487,488\}\]
	\[\min|\hat{\lambda}_{H}^{p,q}| =  |10.166|,\quad \text{index}=\{852,853\}\]

Numerical experiments show that with (\ref{eq:GeneralFormReD}) and (\ref{eq:optstcl_k80}) as a re-discretization stencil for the coarse grid, the solver converges but does not show wavenumber-independent convergence. The iterations required to solve the constant-wavenumber model problem with a point source in the center and Dirichlet boundary conditions are in Table \ref{tab:ReD-optstcl-vs-ReD-O4}. It is comparable to the coarse-grid re-discretization using a standard fourth-order finite scheme (ReD-O4).
\begin{table}[htbp]
    \centering
    \caption{The outer iterations required to solve the constant-wavenumber model problem with a point source in the center and Dirichlet boundary conditions. }
    \label{tab:ReD-optstcl-vs-ReD-O4}
    \scalebox{0.89}{
        \begin{tabular}{lrlll}
        \hline
            Grid size& k & ReD-9ptO2 & ReD-O2 &ReD-O4   \\ \hline
            $129\times 129$& 40 & 15   & 17  & 14\\
            $257\times 257$& 80 & 20   & 28  & 19\\
            $321\times 321$& 100& 27   & 37  & 22\\
        \hline
        \end{tabular}
    }
\end{table}

\clearpage
\section*{Acknowledgments}
We would like to acknowledge the support of the CSC scholarship (No. 202006230087).

\bibliographystyle{siamplain}
\bibliography{Two-level_Deflation_arXiv}

\begin{thebibliography}{10}

\bibitem{babuska1997pollution}
{\sc I.~M. Babuska and S.~A. Sauter}, {\em Is the pollution effect of the {FEM}
  avoidable for the {H}elmholtz equation considering high wave numbers?}, SIAM
  Journal on numerical analysis, 34 (1997), pp.~2392--2423,
  \url{https://doi.org/10.1137/S0036142994269186}.

\bibitem{bayliss1983iterative}
{\sc A.~Bayliss, C.~I. Goldstein, and E.~Turkel}, {\em An iterative method for
  the {H}elmholtz equation}, Journal of Computational Physics, 49 (1983),
  pp.~443--457, \url{https://doi.org/10.1016/0021-9991(83)90139-0}.

\bibitem{boubendir2012quasi}
{\sc Y.~Boubendir, X.~Antoine, and C.~Geuzaine}, {\em A quasi-optimal
  non-overlapping domain decomposition algorithm for the {H}elmholtz equation},
  Journal of Computational Physics, 231 (2012), pp.~262--280,
  \url{https://doi.org/10.1016/j.jcp.2011.08.007}.

\bibitem{calandra2013improved}
{\sc H.~Calandra, S.~Gratton, X.~Pinel, and X.~Vasseur}, {\em An improved
  two-grid preconditioner for the solution of three-dimensional {H}elmholtz
  problems in heterogeneous media}, Numerical Linear Algebra with Applications,
  20 (2013), pp.~663--688, \url{https://doi.org/10.1002/nla.1860}.

\bibitem{calandra2017geometric}
{\sc H.~Calandra, S.~Gratton, and X.~Vasseur}, {\em A geometric multigrid
  preconditioner for the solution of the {H}elmholtz equation in
  three-dimensional heterogeneous media on massively parallel computers}, in
  Modern Solvers for {H}elmholtz Problems, Springer, 2017, pp.~141--155,
  \url{https://doi.org/10.1007/978-3-319-28832-1_6}.

\bibitem{jchen2D2022}
{\sc J.~Chen, V.~Dwarka, and C.~Vuik}, {\em Matrix-free parallel preconditioned
  iterative solvers for the 2{D} helmholtz equation discretized with finite
  differences}, in Scientific Computing in Electrical Engineering, Springer
  International Publishing, 2023.

\bibitem{chen2013source}
{\sc Z.~Chen and X.~Xiang}, {\em A source transfer domain decomposition method
  for {H}elmholtz equations in unbounded domain}, SIAM Journal on Numerical
  Analysis, 51 (2013), pp.~2331--2356, \url{https://doi.org/10.1137/130917144}.

\bibitem{collino2000domain}
{\sc F.~Collino, S.~Ghanemi, and P.~Joly}, {\em Domain decomposition method for
  harmonic wave propagation: a general presentation}, Computer Methods in
  Applied Mechanics and Engineering, 184 (2000), pp.~171--211,
  \url{https://doi.org/10.1016/S0045-7825(99)00228-5}.

\bibitem{DHPC2022}
{\sc {D}elft {H}igh {P}erformance {C}omputing~{C}entre ({DHPC})}, {\em
  {D}elft{B}lue {S}upercomputer ({P}hase 1)}.
\newblock \url{https://www.tudelft.nl/dhpc/ark:/44463/DelftBluePhase1}, 2022.

\bibitem{dostal1988conjugate}
{\sc Z.~Dost{\'a}l}, {\em Conjugate gradient method with preconditioning by
  projector}, International Journal of Computer Mathematics, 23 (1988),
  pp.~315--323.

\bibitem{douglas1998second}
{\sc J.~Douglas~Jr and D.~B. Meade}, {\em Second-order transmission conditions
  for the {H}elmholtz equation}, in Ninth International Conference on Domain
  Decomposition Methods, Citeseer, 1998, pp.~434--440.

\bibitem{DWARKA2022111327}
{\sc V.~Dwarka and C.~Vuik}, {\em Scalable multi-level deflation
  preconditioning for highly indefinite time-harmonic waves}, Journal of
  Computational Physics, 469 (2022), p.~111327,
  \url{https://doi.org/10.1016/j.jcp.2022.111327}.

\bibitem{dwarka2020scalable}
{\sc V.~N. S.~R. Dwarka and C.~Vuik}, {\em Scalable convergence using two-level
  deflation preconditioning for the {H}elmholtz equation}, SIAM Journal on
  Scientific Computing, 42 (2020), pp.~A901--A928,
  \url{https://doi.org/10.1137/18M1192093}.

\bibitem{Eisenstat1983}
{\sc S.~C. Eisenstat, H.~C. Elman, and M.~H. Schultz}, {\em Variational
  iterative methods for nonsymmetric systems of linear equations}, SIAM Journal
  on Numerical Analysis, 20 (1983), pp.~345--357,
  \url{https://doi.org/10.1137/0720023}.

\bibitem{engquist2011sweeping}
{\sc B.~Engquist and L.~Ying}, {\em Sweeping preconditioner for the {H}elmholtz
  equation: moving perfectly matched layers}, Multiscale Modeling \&
  Simulation, 9 (2011), pp.~686--710, \url{https://doi.org/10.1137/100804644}.

\bibitem{erlangga2012iterative}
{\sc Y.~Erlangga and E.~Turkel}, {\em Iterative schemes for high order compact
  discretizations to the exterior {H}elmholtz equation}, ESAIM: Mathematical
  Modelling and Numerical Analysis, 46 (2012), pp.~647--660.

\bibitem{erlangga2008multilevel}
{\sc Y.~A. Erlangga and R.~Nabben}, {\em Multilevel projection-based nested
  {K}rylov iteration for boundary value problems}, SIAM Journal on Scientific
  Computing, 30 (2008), pp.~1572--1595.

\bibitem{erlangga2008MLKM}
{\sc Y.~A. Erlangga and R.~Nabben}, {\em On a multilevel krylov method for the
  {H}elmholtz equation preconditioned by shifted laplacian}, Electronic
  Transactions on Numerical Analysis, 31 (2008), p.~3.

\bibitem{erlangga2006novel}
{\sc Y.~A. Erlangga, C.~W. Oosterlee, and C.~Vuik}, {\em A novel multigrid
  based preconditioner for heterogeneous {H}elmholtz problems}, SIAM Journal on
  Scientific Computing, 27 (2006), pp.~1471--1492,
  \url{https://doi.org/10.1137/040615195}.

\bibitem{erlangga2004class}
{\sc Y.~A. Erlangga, C.~Vuik, and C.~W. Oosterlee}, {\em On a class of
  preconditioners for solving the {H}elmholtz equation}, Applied Numerical
  Mathematics, 50 (2004), pp.~409--425,
  \url{https://doi.org/10.1016/j.apnum.2004.01.009}.

\bibitem{gander2015applying}
{\sc M.~J. Gander, I.~G. Graham, and E.~A. Spence}, {\em Applying {GMRES} to
  the {H}elmholtz equation with shifted laplacian preconditioning: what is the
  largest shift for which wavenumber-independent convergence is guaranteed?},
  Numerische Mathematik, 131 (2015), pp.~567--614,
  \url{https://doi.org/10.1007/s00211-015-0700-2}.

\bibitem{gander2002optimized}
{\sc M.~J. Gander, F.~Magoules, and F.~Nataf}, {\em Optimized {S}chwarz methods
  without overlap for the {H}elmholtz equation}, SIAM Journal on Scientific
  Computing, 24 (2002), pp.~38--60,
  \url{https://doi.org/10.1137/S1064827501387012}.

\bibitem{gander2019class}
{\sc M.~J. Gander and H.~Zhang}, {\em A class of iterative solvers for the
  {H}elmholtz equation: Factorizations, sweeping preconditioners, source
  transfer, single layer potentials, polarized traces, and optimized {S}chwarz
  methods}, Siam Review, 61 (2019), pp.~3--76,
  \url{https://doi.org/10.1137/16M109781X}.

\bibitem{Gordon2013Robust}
{\sc D.~Gordon and R.~Gordon}, {\em Robust and highly scalable parallel
  solution of the {H}elmholtz equation with large wave numbers}, Journal of
  Computational and Applied Mathematics, 237 (2013), pp.~182--196,
  \url{https://doi.org/10.1016/j.cam.2012.07.024}.

\bibitem{kolotilina1998twofold}
{\sc L.~Y. Kolotilina}, {\em Twofold deflation preconditioning of linear
  algebraic systems. i. theory}, Journal of Mathematical Sciences, 89 (1998),
  pp.~1652--1689.

\bibitem{Kononov2007Numerical}
{\sc A.~V. Kononov, C.~D. Riyanti, S.~W. de~Leeuw, C.~W. Oosterlee, and
  C.~Vuik}, {\em Numerical performance of a parallel solution method for a
  heterogeneous 2{D} {H}elmholtz equation}, Computing and Visualization in
  Science, 11 (2007), pp.~139--146,
  \url{https://doi.org/10.1007/s00791-007-0069-6}.

\bibitem{maclachlan2008algebraic}
{\sc S.~P. MacLachlan and C.~W. Oosterlee}, {\em Algebraic multigrid solvers
  for complex-valued matrices}, SIAM Journal on Scientific computing, 30
  (2008), pp.~1548--1571, \url{https://doi.org/10.1137/070687232}.

\bibitem{mansfield1991damped}
{\sc L.~Mansfield}, {\em Damped jacobi preconditioning and coarse grid
  deflation for conjugate gradient iteration on parallel computers}, SIAM
  Journal on Scientific and Statistical Computing, 12 (1991), pp.~1314--1323.

\bibitem{mcinnes1998additive}
{\sc L.~C. McInnes, R.~F. Susan-Resiga, and D.~E. Keyes}, {\em Additive
  {S}chwarz methods with nonreflecting boundary conditions for the parallel
  computation of {H}elmholtz}, Domain Decomposition Methods 10, 218 (1998),
  p.~325, \url{https://doi.org/10.1090/conm/218/03025}.

\bibitem{nicolaides1987deflation}
{\sc R.~A. Nicolaides}, {\em Deflation of conjugate gradients with applications
  to boundary value problems}, SIAM Journal on Numerical Analysis, 24 (1987),
  pp.~355--365.

\bibitem{plessix2003separation}
{\sc R.~E. Plessix and W.~A. Mulder}, {\em Separation-of-variables as a
  preconditioner for an iterative {H}elmholtz solver}, Applied Numerical
  Mathematics, 44 (2003), pp.~385--400,
  \url{https://doi.org/10.1016/S0168-9274(02)00165-4}.

\bibitem{riyanti2007parallel}
{\sc C.~Riyanti, A.~Kononov, Y.~Erlangga, C.~Vuik, C.~Oosterlee, R.-E. Plessix,
  and W.~Mulder}, {\em A parallel multigrid-based preconditioner for the 3{D}
  heterogeneous high-frequency {H}elmholtz equation}, Journal of Computational
  Physics, 224 (2007), pp.~431--448,
  \url{https://doi.org/10.1016/j.jcp.2007.03.033}.

\bibitem{saad1986gmres}
{\sc Y.~Saad and M.~H. Schultz}, {\em {GMRES}: A generalized minimal residual
  algorithm for solving nonsymmetric linear systems}, SIAM Journal on
  Scientific and Statistical Computing, 7 (1986), pp.~856--869,
  \url{https://doi.org/10.1137/0907058}.

\bibitem{saad2000deflated}
{\sc Y.~Saad, M.~Yeung, J.~Erhel, and F.~Guyomarc'h}, {\em A deflated version
  of the conjugate gradient algorithm}, SIAM Journal on Scientific Computing,
  21 (2000), pp.~1909--1926.

\bibitem{schadle2007additive}
{\sc A.~Sch{\"a}dle and L.~Zschiedrich}, {\em Additive {S}chwarz method for
  scattering problems using the {PML} method at interfaces}, in Domain
  Decomposition Methods in Science and Engineering XVI, Springer, 2007,
  pp.~205--212, \url{https://doi.org/10.1007/978-3-540-34469-8_21}.

\bibitem{PhDSheikh2014}
{\sc A.~H. Sheikh}, {\em Development of the {H}elmholtz Solver based on a
  Shifted {L}aplace Preconditioner and a Multigrid Deflation technique},
  thesis, TU Delft, 2014.

\bibitem{sheikh2013convergence}
{\sc A.~H. Sheikh, D.~Lahaye, and C.~Vuik}, {\em On the convergence of shifted
  laplace preconditioner combined with multilevel deflation}, Numerical Linear
  Algebra with Applications, 20 (2013), pp.~645--662.

\bibitem{singer1998high}
{\sc I.~Singer and E.~Turkel}, {\em High-order finite difference methods for
  the {H}elmholtz equation}, Computer methods in applied mechanics and
  engineering, 163 (1998), pp.~343--358.

\bibitem{stolk2013rapidly}
{\sc C.~C. Stolk}, {\em A rapidly converging domain decomposition method for
  the {H}elmholtz equation}, Journal of Computational Physics, 241 (2013),
  pp.~240--252, \url{https://doi.org/10.1016/j.jcp.2013.01.039}.

\bibitem{PhDTang2008}
{\sc J.~M. Tang}, {\em Two-Level Preconditioned Conjugate Gradient Methods with
  Applications to Bubbly Flow Problem}, thesis, TU Delft, 2008.

\bibitem{tang2008two}
{\sc J.~M. Tang}, {\em Two-level preconditioned conjugate gradient methods with
  applications to bubbly flow problems},  (2008).

\bibitem{taus2020sweeps}
{\sc M.~Taus, L.~Zepeda-N{\'u}{\~n}ez, R.~J. Hewett, and L.~Demanet}, {\em
  L-sweeps: A scalable, parallel preconditioner for the high-frequency
  {H}elmholtz equation}, Journal of Computational Physics, 420 (2020),
  p.~109706, \url{https://doi.org/10.1016/j.jcp.2020.109706}.

\bibitem{toselli1999overlapping}
{\sc A.~Toselli}, {\em Overlapping methods with perfectly matched layers for
  the solution of the {H}elmholtz equation}, in Eleventh International
  Conference on Domain Decomposition Methods, C. Lai, P. Bjorstad, M. Cross,
  and O. Widlund, eds., DDM. org, Citeseer, 1999, pp.~551--558.

\bibitem{Gijzen2007}
{\sc M.~B. van Gijzen, Y.~A. Erlangga, and C.~Vuik}, {\em Spectral analysis of
  the discrete {H}elmholtz operator preconditioned with a shifted laplacian},
  SIAM Journal on Scientific Computing, 29 (2007), pp.~1942--1958,
  \url{https://doi.org/10.1137/060661491}.

\bibitem{Versteeg_1991_ME}
{\em The {Marmousi} experience}, Proceedings of the 1990 EAEG workshop on
  Practical Aspects of Seismic Data Inversion, Eur. Ass. Expl. Geophys, 1991.

\end{thebibliography}
\end{document}